\documentclass[12pt]{article}
\usepackage{amssymb,latexsym, amsmath}
\usepackage[francais]{babel}

\def\cqfd{%
\mbox{ }%
\nolinebreak%
\hfill%
\rule{2.5mm}{2.5mm}%
\medbreak%
\par%
}

\def\N{\mathbb{N}}

\def\R{\mathbb{R}}
\def\C{\mathbb{C}}

\newtheorem{theo}{\quad \  \textsc{Th\'eor\`eme}}[section]
\newtheorem{pro}[theo]{\quad \ \textsc{Proposition}}

\newtheorem{lem}[theo]{\quad \ \textsc{Lemme}}

\title{Fonctions maximales centr\'ees de Hardy-Littlewood pour les op\'erateurs de Grushin}
\author{Hong-Quan LI}
\date{}
\begin{document}
\renewcommand{\theequation}{\thesection.\arabic{equation}}
\setcounter{equation}{0} \maketitle

\vspace{-1.0cm}

\bigskip

{\bf R\'esum\'e.} Consid\'erons l'op\'erateur de Grushin sur $\R_x^n
\times \R_u$, $\Delta_G = \sum_{i = 1}^n \frac{\partial^2}{\partial
x_i^2} + (\sum_{i = 1}^n x_i^2) \frac{\partial^2}{\partial u^2}$.
Soient $d_{CC}$ la distance de Carnot-Carath\'eodory associ\'ee,
$d_K$ une pseudo-distance li\'ee \`a la solution fondamentale de
$\Delta_G$. On montre qu'il existe une constante $A
> 0$, ind\'ependante de $n$,  telle que pour toute $f \in
L^1(\R_x^n \times \R_u, dxdu)$, on a $\| M f \|_{L^{1, \infty}} \leq
A n \| f \|_1$, o\`u $M$ d\'esigne la fonction maximale centr\'ee de
Hardy-Littlewood d\'efinie soit par $d_{CC}$, soit par $d_K$. On
trouve une relation \'etroite entre ce sujet et la fonction de
Green.

\medskip

{\bf Mathematics Subject Classification (2000) :} {\bf 42B25,
43A80}

\medskip

{\bf Key words and phrases :} Centered Hardy-Littlewood maximal
function; Grushin operators; Poisson kernel; Carnot-Carath\'eodory
distance; Green function

\medskip

\renewcommand{\theequation}{\thesection.\arabic{equation}}
\section{Introduction}
\setcounter{equation}{0}
\medskip

Consid\'erons l'op\'erateur de Grushin sur $\R_x^n \times \R_u$, qui
a \'et\'e initialement \'etudi\'e par Grushin (voir par exemple
\cite{G71})
\begin{eqnarray*}
\Delta_G = \sum_{i = 1}^n \frac{\partial^2}{\partial x_i^2} + |x|^2
\frac{\partial^2}{\partial u^2} = \sum_{i = 1}^n (X_i^2 + U_i^2),
\end{eqnarray*}
avec
\begin{eqnarray*}
|x|^2 = \sum_{i = 1}^n x_i^2, \quad X_i = \frac{\partial}{\partial
x_i}, \quad U_i = x_i \frac{\partial}{\partial u},
\end{eqnarray*}
o\`u les champs de vecteurs infiniment diff\'erentiables $\{ X_i,
U_i \}_{1 \leq i \leq n}$ satisfont la condition de H\"ormander.

Notons $d_{CC}$ la distance de Carnot-Carath\'eodory associ\'ee \`a
$\{X_1, \cdots, X_n, U_1, \cdots, U_n \}$, voir par exemple
\cite{VSCC92}. Pour simplifier les notations, on pose
\begin{eqnarray} \label{NNN1}
g = (x, u), g' = (x', u') \in \R^n_x \times \R_u, \quad s = | u - u'|, \nonumber \\
R^2 = |x|^2 + |x'|^2, \quad a = \frac{2 x \cdot x'}{R^2} \in [-1,
1],
\end{eqnarray}
avec $x \cdot x'$ le produit scalaire sur $\R_x^n$; d\'efinissons
pour $-\pi < \varphi < \pi$,
\begin{eqnarray} \label{NNN7}
\mu(a; \varphi) = \frac{\varphi}{\sin^2{\varphi}} - \cot{\varphi} +
a \frac{1 - \varphi \cot{\varphi}}{\sin{\varphi}}.
\end{eqnarray}

On sait que $\mu(a; \cdot)$ est strictement croissante avec $\mu(a;
0) = 0$ et est un diff\'eomophisme de $]-\pi, \pi[$ sur $\R$ pour
$-1 < a \leq 1$, et est un diff\'eomophisme de $]-\pi, \pi[$ sur
$]-\frac{\pi}{2}, \frac{\pi}{2}[$ pour $a = -1$, voir les Lemmas
2.2-2.5 de \cite{P07}. Notons $\mu^{-1}(a; \cdot)$ sa fonction
r\'eciproque.
On a (voir le Theorem 2.6 de \cite{P07})  \\
(1) dans le cas o\`u $x = -x'$ avec $2 s \geq \pi |x|^2$, $d_{CC}(g,
g') = \sqrt{2 \pi s}$. \\
(2) dans d'autre cas,
\begin{eqnarray} \label{FDCC}
d^2_{CC}(g, g') = \Big( \frac{\theta}{\sin{\theta}} \Big)^2 R^2 (1 -
a \cos{\theta}) \ \mbox{avec} \ \theta = \mu^{-1}\Big(a; \frac{2
s}{R^2}\Big).
\end{eqnarray}

Lorsque $g = 0$, on retrouve la distance de Carnot-Carath\'eodory
sur le groupe de Heisenberg (voir \cite{BGG00}).

On d\'efinit maintenant une pseudo-distance $d_K$ li\'ee \`a la
solution fondamentale de $\Delta_G$ (c'est-\`a-dire la fonction de
Green, voir la section \S \ref{S2}), qui est cruciale pour cet
article :
\begin{eqnarray} \label{FDK}
d_K(g, g') = \Big( \sqrt{R^4 + (2 s)^2} - 2 x \cdot x'
\Big)^{\frac{1}{2}}.
\end{eqnarray}

Lorsque $g = 0$, on retrouve la norme de Kor\'anyi sur le groupe de
Heisenberg (voir par exemple \cite{F73}).

On voit que $d_K$ et $d_{CC}$ sont \'equivalentes, plus
pr\'ecisement, il existe une constante $C > 1$, ind\'ependante de
$n$, telle que
\begin{eqnarray} \label{EKCC}
d_K(g, g') \leq d_{CC}(g, g') \leq C d_K(g, g'), \qquad \forall g,
g'.
\end{eqnarray}
En effet, l'in\'egalit\'e \`a gauche est donn\'ee par la Proposition
\ref{LL1} de cet article. Pour v\'erifier l'in\'egalit\'e \`a
droite, il suffit d'utiliser la Proposition 5.1 de \cite{RS08} et
d'observer que les (pseudo-)distances sont ind\'ependante de
dimension. Ce sera tr\`es int\'eressant de savoir si $d_K$ est une
distance sur $\R^n \times \R$.

Dans la suite, on note $B_K(g, r)$ la boule ouverte de centre $g$ et
de rayon $r > 0$ d\'efinie par la pseudo-distance $d_K$, et
$B_{CC}(g, r)$ celle d\'efinie par la distance de
Carnot-Carath\'eodory. Si $E$ est un ensemble mesurable, alors on
note $|E|$ son volume et $\chi_E$ sa fonction caract\'eristique. On
montrera que pour tout $g = (x, u)$ et tout $r > 0$, on a
\begin{eqnarray} \label{EVK}
\frac{1}{8} r^{n+1} (r + |x|) B(\frac{n}{2}, \frac{3}{2}) | S^{n-1}|
\leq |B_K(g, r)| \leq r^{n+1} (r + |x|) B(\frac{n}{2}, \frac{3}{2})
| S^{n-1}|,
\end{eqnarray}
o\`u $|\mathrm{S}^{n -1}|$ d\'esigne la surface de la sph\`ere
unitaire dans $\R^{n}$, $B(\frac{n}{2}, \frac{3}{2})$ \'etant la
fonction de Beta d'indices $\frac{n}{2}$ et $\frac{3}{2}$. Il existe
une structure de dilatation sur $(\R^n \times \R, d, dg)$ avec $d =
d_K$ ou bien $d_{CC}$ : \\
Pour tout $r > 0$, d\'efinissons
\begin{eqnarray*}
\delta_r : & \R^n \times \R \longrightarrow \R^n \times \R \\
\mbox{} & \delta_r(x, u) = (r x, r^2 u);
\end{eqnarray*}
et pour $A \subset \R^n \times \R$, d\'efinissons $\delta_r A =
\{\delta_r(g); g \in A \}$. Alors, pour tous $g, g'$ et tout $r >
0$, on a
\begin{eqnarray*}
d(\delta_r(g), \delta_r(g')) = r d(g, g'), \quad B(g, r) = \delta_r
B(\delta_{r^{-1}}g, 1), \quad |B(g, r)| = r^{n+2}
|B(\delta_{r^{-1}}g, 1)|,
\end{eqnarray*}
avec $B = B_K$ ou bien $B = B_{CC}$. Cependant, $(\R^n \times \R, d,
dg)$ n'a pas de la propi\'et\'e d'invariance par translation.

Par les Propositions \ref{LL1} et \ref{LL2} ci-apr\`es, il existe
une constante $c > 0$ telle que pour tout $n
\in \N^*$, tout $g \in \R^n \times \R$ et tout $r > 0$, on a
\begin{eqnarray} \label{EVKCC}
c |B_{K}(g, r)| \leq |B_{CC}(g, r)| \leq |B_{K}(g, r)|.
\end{eqnarray}

En particulier, \eqref{EVK} et \eqref{EVKCC} nous disent que $(\R^n
\times \R, d, dg)$ est un espace de nature homog\`ene au sens de
Coifman-Weiss (voir \cite{CW71}).

Pour $f \in L_{loc}^1(\R^n \times \R)$, on d\'efinit les deux
fonctions maximales centr\'ees de Hardy-Littlewood $M_Kf$ et
$M_{CC}f$ respectivement par
\begin{eqnarray*}
M_Kf(g) &=& \sup_{r > 0} |B_K(g, r)|^{-1} \int_{B_K(g, r)} |f(g')|
\, dg', \quad  \forall g, \\
M_{CC}f(g) &=& \sup_{r > 0} |B_{CC}(g, r)|^{-1} \int_{B_{CC}(g, r)}
|f(g')| \, dg', \quad \forall g.
\end{eqnarray*}

Notre r\'esultat principal est le

\begin{theo} \label{TH}
Il existe une constante $L > 0$ telle que pour tout $n \in \N^*$, on
a
\begin{eqnarray} \label{E1N}
\| M f \|_{L^{1, \infty}} \leq L n \| f \|_1, \quad \forall f \in
L^1(\R^n \times \R),
\end{eqnarray}
avec $M = M_K$ ou bien $M = M_{CC}$.
\end{theo}

Ce type d'estimation a \'et\'e obtenu par Stein et Str\"omberg dans
le cadre des espaces euclidiens pour la fonction maximale standard
de Hardy-Littlewood (voir \cite{SS83}). Remarquons que l'estimation
pr\'ec\'edente a \'et\'e obtenue dans \cite{Li09} sur les groupes de
Heisenberg pour la fonction maximale d\'efinie par la distance de
Carnot-Carath\'eodory ou bien par celle de Kor\'anyi.

Par \eqref{FDK}, \eqref{EVK} et \eqref{EVKCC}, on peut montrer que
$(\R^n \times \R, d_{CC}, dg)$ satisfait la propri\'et\'e de
``strong $n+1$-microdoubling with constant $L_1$'' au sens de
\cite{NT09}, i.e.
\begin{eqnarray*}
|B_{CC}(g', (1 + \frac{1}{n+1}) r)| \leq L_1 |B_{CC}(g, r)|, \quad
\forall g, r > 0, g' \in B_{CC}(g, r).
\end{eqnarray*}

Par le Corollary 1.2 de \cite{NT09}, qui est motiv\'ee par
\cite{SS83} et valable dans une situation tr\`es g\'en\'erale, il
existe une constante $C(L_1)
> 0$, ind\'ependente de $n$, telle que
\begin{eqnarray} \label{E2}
\| M_{CC} f \|_{L^{1, \infty}} \leq C(L_1) (n+1) \ln{(n+1)} \| f
\|_{L^1}, \quad \forall f \in L^1.
\end{eqnarray}

Pour d'autres travaux au sujet des estimations de type
\eqref{E2} et \eqref{E1N}, voir \cite{SS83}, \cite{NT09}
et leurs r\'ef\'erences.

L'id\'ee principale de la d\'emonstration du Th\'eor\`eme \ref{TH}
est d'utiliser ``the Hopf-Dunford-Schwartz maximal ergodic theorem''
(voir \cite{DS58} pp.690-691) comme dans \cite{SS83}. Afin de
comprendre d'o\`u vient l'estimation \eqref{E1N}, on explique
bri\`evement, au point de vue de la fonction de Green, la preuve de
Stein-Str\"omberg dans le cadre des espaces euclidiens : \\
en notant $e^{h \Delta_{\R^n}}$ ($h > 0$) le semi-groupe de la
chaleur, par ``the Hopf-Dunford-Schwartz maximal ergodic theorem''
et la structure de dilatation sur $\R^n$, il suffit de montrer qu'il
existe une constante $A > 0$, ind\'ependante de $n$ ($\geq 3$, si
l'on veut), telle que pour un certain $s(n) > 0$, tout $x \in \R^n$
et pour tout $x \neq x' \in B_{\R^n}(x, 1)$, on a
\begin{eqnarray*}
|B_{\R^n}(x, 1)|^{-1} &\leq& A n \frac{1}{s(n)} \int_0^{s(n)} e^{h
\Delta_{\R^n}}(x, x') \, dh \\
&=& \frac{A n}{s(n)} \Big[ (-\Delta_{\R^n})^{-1} -
(-\Delta_{\R^n})^{-1} e^{s(n) \Delta_{\R^n}} \Big](x, x'),
\end{eqnarray*}
par la formule explicite du noyau de la chaleur sur $\R^n$, en
choisissant $s(n) = \frac{1}{n}$, on voit
qu'il existe une constante $c > 0$, ind\'ependante de $n$, telle que
\begin{eqnarray} \label{EGGG1}
\Big[ (-\Delta_{\R^n})^{-1} - (-\Delta_{\R^n})^{-1} e^{s(n)
\Delta_{\R^n}} \Big](x, x') \geq c (-\Delta_{\R^n})^{-1}(x, x'),
 x \neq x' \in B_{\R^n}(x, 1).
\end{eqnarray}
Pour terminer la preuve, il suffit d'utiliser le fait que pour $n
\geq 3$,
\begin{eqnarray*}
(-\Delta_{\R^n})^{-1}(x, x') = \frac{|x - x'|^{2 - n}}{n (n - 2)}
|B_{\R^n}(x, 1)|^{-1}.
\end{eqnarray*}
Remarquons que $(-\Delta_{\R^n})^{-1} e^{s(n) \Delta_{\R^n}}$ est
positif, on a
\begin{eqnarray*}
\Big[ (-\Delta_{\R^n})^{-1} - (-\Delta_{\R^n})^{-1} e^{s(n)
\Delta_{\R^n}} \Big](x, x') \leq (-\Delta_{\R^n})^{-1}(x, x').
\end{eqnarray*}

En gros, il faut prendre $\phi(n) = n$ telle que
\begin{eqnarray} \label{GVEE}
\inf_{x \neq x' \in B_{\R^n}(x, 1), n \geq 3} \phi(n) |B_{\R^n}(x,
1)| n (-\Delta_{\R^n})^{-1}(x, x') > 0,
\end{eqnarray}
et on obtient \eqref{E1N} avec $n = \phi(n)$. On remarque que le
r\'esultat de \cite{Li09} peut aussi \^etre expliqu\'e par une
estimation de type \eqref{GVEE}. Naturellement, on s'int\'eresse \`a
\'etudier la fonction de Green pour l'op\'erateur de Grushin et \`a
\'etudier une estimation de type \eqref{GVEE}. \c{C}a nous explique
pourquoi la pseudo-distance s'introduit et pourquoi on a le
th\'eor\`eme \ref{TH}. Evidemment, cette explication n'est pas
s\'erieuse puisqu'on ne sait pas comment \'etablir l'estimation de
type \eqref{EGGG1} dans le cadre des groupes de Heisenberg et des
op\'erateurs de Grushin. Remarquons qu'on a obtenu dans \cite{Li07}
et dans \cite{Li10} des estimations asymptotiques du noyau de la
chaleur sur les groupes de Heisenberg et pour $\Delta_G$
respectivement; cependant, ces estimations ne sont pas
uniform\'ement en dimension (ce sera tr\`es int\'eressant de
consid\'erer ce probl\`eme et d'\'etudier l'estimation de type
\eqref{EGGG1}). Donc, on utilisera le noyau de Poisson comme dans
\cite{Li09} au lieu du noyau de la chaleur utilis\'e dans
\cite{SS83}. Dans les grandes lignes, il s'agit d'une estimation de
type
\begin{align} \label{ITI}
\inf_{x \neq x' \in B(x, 1), n \geq 3} \phi(n) |B(x, 1)| \sqrt{n}
(-\Delta)^{-\frac{1}{2}}(x, x') > 0,
\end{align}
o\`u le choix de $\sqrt{n}$ est optimal \`a part d'une constante
universelle comme on a expliqu\'e dans \cite{Li09}. Voir l'appendice
pour plus d'explication.

Dans \cite{LL10}, on s'adapte les m\'ethodes de cet article et de
\cite{Li09} pour qu'ils soient valables dans des situations beaucoup
plus compliqu\'ees. Plus pr\'ecisement, sur les espaces
hyperboliques r\'eel de dimension $n$ ($n \geq 2$), $\mathbb{H}^n$,
qui sont \`a croissance exponentielle du volume, on montre qu'il
existe une constante $L > 0$, ind\'ependante de $n$, telle que
\begin{align*}
\| M \|_{L^1 \longrightarrow L^{1, \infty}} \leq L n \ln{n},
\end{align*}
o\`u $M$ d\'esigne la fonction maximale centr\'ee de
Hardy-Littlewood sur $\mathbb{H}^n$.

Cet article est organis\'e de la fa\c con suivante : on consid\'era
dans la section \ref{S2} la fonction de Green pour l'op\'erateur de
Grushin et la pseudo-distance $d_K$. On \'etudiera les estimations
uniform\'ement asymptotiques du noyau de Poisson dans la section
\ref{S3}. La d\'emonstration du Th\'eor\`eme \ref{TH} sera donn\'ee
dans la section \ref{S4} pour $M = M_K$ et dans la section \ref{S5}
pour $M = M_{CC}$.

\subsection {Notations}

Dans toute la suite, $c$, $A$, etc. d\'esigneront des constantes
universelles qui peuvent changer d'une ligne
\`a l'autre.

Pour deux fonctions $f$ et $g$, on dit que $f = O(g)$ s'il existe
une constante $c > 0$ telle que $|f| \leq c |g|$; que $f = o(g)$ si
$\lim \frac{f}{g} = 0$; et que $f \sim g$ s'il existe une constante
$A > 1$ telle que $A^{-1} f \leq g \leq A f$.

\medskip

\renewcommand{\theequation}{\thesection.\arabic{equation}}
\section{La fonction de Green pour l'op\'erateur de Grushin}\label{S2}
\setcounter{equation}{0}

\medskip

Apr\`es avoir fini ce travail, on trouve que dans \cite{BGG99},
Beals et al. ont obtenu une expression pour la fonction de Green
dans une situation un peu plus g\'en\'erale. Cependant, leur
expression n'est ni claire, ni pratique pour cet article, surtout
elle n'offert aucune information pour s'introduire naturellement la
pseudo-distance $d_K$, qui est essentiel pour ce travail. On
utilisera une autre m\'ethode pour donner la fonction de Green.

Soit $p_h = p_h^{(n)}$ ($h > 0$) le noyau de la chaleur
(c'est-\`a-dire le noyau int\'egral de $e^{h \Delta_G}$) sur $(\R^n
\times \R, \Delta_G, dg)$. L'expression explicite suivante de $p_h$
a \'et\'e obtenue par Paulat dans \cite{P07} :
\begin{eqnarray} \label{ehk1}
p_h((x, u), (x', u')) = (4 \pi h)^{-\frac{n}{2} - 1} K(\frac{1}{4 h}
(|x|^2 + |x'|^2), \frac{2 x \cdot x'}{|x|^2 + |x'|^2}; \frac{1}{4 h}
| u - u' |),
\end{eqnarray}
o\`u la fonction $K$ est d\'efinie sur $\R_+ \times [-1, 1] \times
\R_+$ par
\begin{eqnarray} \label{ehkH1}
K(s_1, s_2; s_3) = \int_{\R} \Big( \frac{\lambda}{\sinh{\lambda}}
\Big)^{\frac{n}{2}} \exp{\Big( 2 i \lambda s_3 - s_1 ( \lambda
\coth{\lambda} - \frac{\lambda}{\sinh{\lambda}} s_2) \Big)} \,
d\lambda.
\end{eqnarray}

La fonction de Green est donn\'ee par
\begin{eqnarray*}
(-\Delta_G)^{-1}(g, g') = \int_0^{+\infty} p_h(g, g') \, dh;
\end{eqnarray*}
notons
\begin{eqnarray*}
X = X(g, g', \lambda) = \int_0^{+\infty} h^{-\frac{n}{2} - 1}
\exp\Big\{ \frac{1}{4 h} \Big( i (2 s \lambda) - R^2 (\lambda
\coth{\lambda} - \frac{\lambda}{\sinh{\lambda}} a) \Big) \Big\} \, dh,
\end{eqnarray*}
on a (formulement)
\begin{eqnarray*}
(-\Delta_G)^{-1}(g, g') = (4 \pi)^{-\frac{n}{2}-1} \int_{\R} \Big(
\frac{\lambda}{\sinh{\lambda}} \Big)^{\frac{n}{2}} X \, d\lambda.
\end{eqnarray*}

Par le changement de variable
\begin{eqnarray*}
\gamma = \frac{R^2 (\lambda \coth{\lambda} -
\frac{\lambda}{\sinh{\lambda}} a) - i (2 s \lambda)}{4 h},
\end{eqnarray*}
on voit que
\begin{eqnarray*}
X &=& \Big[ \frac{R^2 (\lambda \coth{\lambda} -
\frac{\lambda}{\sinh{\lambda}} a) - i (2 s \lambda)}{4}
\Big]^{-\frac{n}{2}} \int_0^{+\infty} \gamma^{\frac{n}{2} - 1} e^{-\gamma} \, d\gamma \\
&=& \Big[ \frac{4}{R^2 (\lambda \coth{\lambda} -
\frac{\lambda}{\sinh{\lambda}} a) - i (2 s \lambda)} \Big]^{
\frac{n}{2}} \Gamma(\frac{n}{2}),
\end{eqnarray*}
donc,
\begin{eqnarray*}
(-\Delta_G)^{-1}(g, g') = \frac{\Gamma(\frac{n}{2})}{4 \pi^{
\frac{n}{2} + 1}} \int_{\R} \Big[ R^2 \cosh{\lambda} - i (2 s
\sinh{\lambda}) - R^2 a \Big]^{-\frac{n}{2}} \, d\lambda.
\end{eqnarray*}

Posons
\begin{eqnarray*}
D_K(R^2, 2s) = \Big[ R^4 + (2 s)^2 \Big]^{\frac{1}{4}},
\end{eqnarray*}
et $0 \leq \phi \leq \frac{\pi}{2}$ telle que
\begin{eqnarray*}
e^{- i \phi} = D_K^{-2}(R^2, 2s) [R^2 - i (2s)].
\end{eqnarray*}

On a alors
\begin{eqnarray} \label{EGF1}
R^2 \cosh{\lambda} - i (2 s \sinh{\lambda}) = D_K^{2}(R^2, 2s)
\cosh{(\lambda - i \phi)},
\end{eqnarray}
et formulement
\begin{eqnarray*}
(-\Delta_G)^{-1}(g, g') &=& \frac{\Gamma(\frac{n}{2})}{4 \pi^{
\frac{n}{2} + 1}} \int_{\R} \Big[ D_K^{2}(R^2, 2s) \cosh{(\lambda -
i \phi)} - R^2 a \Big]^{-\frac{n}{2}} \, d\lambda \\
&=& \frac{\Gamma(\frac{n}{2})}{4 \pi^{ \frac{n}{2} + 1}} \int_{\R}
\Big[ D_K^{2}(R^2, 2s) \cosh{\lambda} - R^2 a \Big]^{-\frac{n}{2}}
\, d\lambda.
\end{eqnarray*}

En utilisant l'expression int\'egrale de la fonction de Legendre
$P_{-\frac{1}{2}}^{-\frac{n-1}{2}}$ (voir \cite{EMOT53} p.156), on
peut obtenir une expression pour $(-\Delta_G)^{-1}$. Cependant, ce
n'est pas n\'ecessaire pour cet article. Remarquons que la
dimension homog\`ene de $(\R^n \times \R, \Delta_G, dg)$ est de $n +
2$. En rappelant l'estimation classique de la fonction de Green
hypoelliptique (voir par exemple \cite{SC84}, \cite{NSW85}, \cite{FSC86}
et \cite{JSC87}, qui peut aussi \^etre expliqu\'ee facilement par les estimations
classiques du noyau de la chaleur), on \'ecrit
naturellement
\begin{eqnarray*}
(-\Delta_G)^{-1}(g, g') &=& \frac{\Gamma(\frac{n}{2})}{4 \pi^{
\frac{n}{2} + 1}} \int_{\R} \Big[ ( D_K^{2}(R^2, 2s) - R^2 a) + 2
D_K^{2}(R^2, 2s) \sinh^2{\frac{\lambda}{2}} \Big]^{-\frac{n}{2}} \,
d\lambda \\
&=& \frac{\Gamma(\frac{n}{2})}{4 \pi^{ \frac{n}{2} + 1}} d_K^{-n}(g,
g') \int_{\R} \Big[ 1 + \frac{2 D_K^{2}(R^2, 2s)}{d_K^2(g, g')}
\sinh^2{\frac{\lambda}{2}} \Big]^{-\frac{n}{2}} \, d\lambda \\
&=& \frac{\Gamma(\frac{n}{2})}{\pi^{ \frac{n}{2} + 1}} d_K^{-n}(g,
g') \int_0^{+\infty} \Big[ 1 + \frac{2 D_K^{2}(R^2, 2s)}{d_K^2(g,
g')} \sinh^2{\lambda} \Big]^{-\frac{n}{2}} \, d\lambda.
\end{eqnarray*}

Pour simplifier les notations, posons
\begin{eqnarray*}
Y = \int_0^{+\infty} \Big[ 1 + \frac{2 D_K^{2}(R^2, 2s)}{d_K^2(g,
g')} \sinh^2{\lambda} \Big]^{-\frac{n}{2}} \, d\lambda.
\end{eqnarray*}
On voit facilement que
\begin{eqnarray*}
\frac{2 D_K^{2}(R^2, 2s)}{d_K^2(g, g')} \geq 1.
\end{eqnarray*}
On montrera qu'il existe une constante $c
> 1$, ind\'ependante de $(n, \frac{2 D_K^{2}(R^2, 2s)}{d_K^2(g,
g')})$, telle que pour tout $n \geq 2$ et tout $\frac{2 D_K^{2}(R^2,
2s)}{d_K^2(g, g')} \geq 1$,
\begin{eqnarray}
c^{-1} n^{-\frac{1}{2}} \frac{d_K(g, g')}{D_K(R^2, 2s)} \leq Y \leq
c n^{-\frac{1}{2}} \frac{d_K(g, g')}{D_K(R^2, 2s)}.
\end{eqnarray}

En fait, puisque $\sinh{\lambda} \geq \lambda$ pour tout $\lambda
\geq 0$, on a
\begin{eqnarray*}
Y &\leq& \int_0^{+\infty} \Big[ 1 + \frac{2 D_K^{2}(R^2,
2s)}{d_K^2(g, g')} \lambda^2 \Big]^{-\frac{n}{2}} \, d\lambda \\
&=& \frac{1}{2 \sqrt{2}} \frac{d_K(g, g')}{D_K(R^2, 2s)}
B(\frac{n-1}{2}, \frac{1}{2}),
\end{eqnarray*}
voir \cite{EMOT53} p.10 (16). On sait qu'il existe une constante $C
> 0$ telle que $B(\frac{n-1}{2},
\frac{1}{2}) \leq C n^{-\frac{1}{2}}$ pour tout $n \geq 2$. Donc, on
obtient la majoration de $Y$.

D'ailleurs, on a
\begin{eqnarray*}
Y = \int_0^{+\infty} \exp{\Big[ -\frac{n}{2} \ln{(1 + \frac{2
D_K^{2}(R^2, 2s)}{d_K^2(g, g')} \sinh^2{\lambda})} \Big]} \,
d\lambda.
\end{eqnarray*}
Or,
\begin{eqnarray*}
\sinh{\lambda} \leq 2 \lambda, \quad \forall 0 \leq \lambda \leq 1,
\quad \ln{(1 + s)} \leq s, \quad \forall s \geq 0,
\end{eqnarray*}
on a
\begin{eqnarray*}
Y &\geq& \int_0^1 \exp{\Big[ -\frac{n}{2} \ln{(1 + \frac{2
D_K^{2}(R^2, 2s)}{d_K^2(g, g')} 4 \lambda^2)} \Big]} \,
d\lambda \\
&\geq& \int_0^{(4 \frac{2 D_K^{2}(R^2, 2s)}{d_K^2(g,
g')})^{-\frac{1}{2}}} \exp{\Big[ - 4 \frac{n}{2} \frac{2
D_K^{2}(R^2, 2s)}{d_K^2(g, g')} \lambda^2 \Big]} \, d\lambda \\
&\geq& 2^{-1} n^{-\frac{1}{2}} \frac{d_K(g, g')}{D_K(R^2, 2s)}
\int_0^1 e^{-h^2} \, dh.
\end{eqnarray*}

Par cons\'equent, il existe une constante $A_1 > 1$ telle que pour
tout $n \geq 2$ et tous $g \neq g'$, on a
\begin{eqnarray} \label{GEV}
A_1^{-1} \leq \frac{(-\Delta_G)^{-1}(g,
g')}{\frac{\Gamma(\frac{n}{2})}{4 \pi^{ \frac{n}{2} + 1}}
d_K^{-n}(g, g') n^{-\frac{1}{2}} \frac{d_K(g, g')}{D_K(R^2, 2s)}}
\leq A_1.
\end{eqnarray}

\subsection{Estimation de $|B_K(g, 1)|$}

Le but de cette sous-section est de montrer l'estimation
\eqref{EVK}. Par la structure de dilatation, il suffit de
consid\'erer le cas o\`u $r = 1$.

Rappelons que
\begin{eqnarray*}
d_K^2((x, u), (x', u')) = \sqrt{(|x|^2 + |x'|^2)^2 + (2 |u' - u|)^2}
- 2 x' \cdot x,
\end{eqnarray*}
on constate d'abord que
\begin{eqnarray*}
(x', u') \in B_K((x, u), 1) \Longleftrightarrow (2 |u' - u|)^2 < (1
+ 2 x' \cdot x)^2 - (|x|^2 + |x'|^2)^2,
\end{eqnarray*}
puis, par le fait que
\begin{eqnarray*}
(1 + 2 x' \cdot x)^2 - (|x|^2 + |x'|^2)^2 = (1 - |x' - x|^2) \cdot
(1 + |x' + x|^2),
\end{eqnarray*}
que
\begin{eqnarray} \label{DDKK1}
(x', u') \in B_K((x, u), 1) \Longleftrightarrow  |x' - x| < 1, 2 |u' - u| < \sqrt{1 - |x' - x|^2} \cdot \sqrt{1 + |x' + x|^2}.
\end{eqnarray}

Par cons\'equent,
\begin{eqnarray} \label{VE1}
|B_K((x, u), 1)| &=& \int_{|x' - x| < 1} \sqrt{1 - |x' - x|^2} \cdot
\sqrt{1 + |x' + x|^2} \, dx' \nonumber \\
&=& \int_{|z| < 1} \sqrt{1 - |z|^2} \cdot \sqrt{1 + |2 x - z|^2} \,
dz,
\end{eqnarray}
par le changement de variable $x' = x - z$. Comme
\begin{eqnarray*}
\max(2 |x| - 1, 0) \leq |2 x - z| \leq 2 |x| + 1, \quad \forall |z|
< 1,
\end{eqnarray*}
on voit facilement que
\begin{eqnarray} \label{VE2}
\frac{1}{4} (1 + |x|) \leq \sqrt{1 + |2 x - z|^2} \leq 2 (1 + |x|).
\end{eqnarray}

D'ailleurs,
\begin{eqnarray*}
\int_{|z| < 1} \sqrt{1 - |z|^2} \, dz &=& \int_0^1 r^{n-1} \sqrt{1 -
r^2} \, dr |S^{n-1}| \\
&=& \frac{1}{2} |S^{n-1}| \int_0^1 h^{\frac{n}{2} - 1} \sqrt{1 - h}
\, dh = \frac{1}{2} |S^{n-1}| B(\frac{n}{2}, \frac{3}{2}).
\end{eqnarray*}

\eqref{VE1} et \eqref{VE2} impliquent \eqref{EVK}.

Remarquons que
\begin{eqnarray} \label{EDK}
d_K(g, g') < 1 \Longrightarrow D_K(R^2, 2 s) \leq \sqrt{d_K^2(g, g')
+ R^2} \leq 1 + (|x| + |x'|) \leq 2 (1 + |x|),
\end{eqnarray}
par le fait que $|S^{n-1}| = 2
\frac{\pi^{\frac{n}{2}}}{\Gamma(\frac{n}{2})}$, on d\'eduit de
\eqref{EVK} et de \eqref{GEV} qu'une estimation de type \eqref{GVEE}
avec $\phi(n) = n$, i.e.
\begin{eqnarray*}
\inf_{g \neq g' \in B_K(g, 1), n \geq 2} n |B_K(g, 1)| n
(-\Delta_G)^{-1}(g, g') > 0,
\end{eqnarray*}
qui est le point de d\'epart de cet article.

\medskip

\renewcommand{\theequation}{\thesection.\arabic{equation}}
\section{Estimations uniform\'ement asymptotiques du noyau de Poisson pour
$(\R^n \times \R, \Delta_G, dg)$}\label{S3}
\setcounter{equation}{0}

\medskip

On suit la strat\'egie de \cite{Li09}. Il faut dire qu'une
estimation inf\'erieure uniform\'ement du noyau de Poisson est
suffisante pour d\'emontrer le Th\'eor\`eme \ref{TH}.

Soit $P_h = P_h^{(n)}$ ($h > 0$) le noyau de Poisson (c'est-\`a-dire
le noyau int\'egral de $e^{- h \sqrt{-\Delta_G}}$). Par convention,
on note $P = P_1$. Notons $Q = n + 2$. Par la structure de
dilatation, on a
\begin{eqnarray*}
P_h(g, g') &=& h^{-Q} P(\delta_{h^{-1}}(g), \delta_{h^{-1}}(g')),
\quad \forall h > 0, g, g' \in \R^n \times \R.
\end{eqnarray*}

On a les estimations uniform\'ement asymptotiques de $P$, donc
celles de $P_h$, comme suit :

\begin{pro} \label{PP1}
Soit $U \gg 1$. Lorsque $n \longrightarrow + \infty$, pour $d_K(g,
g') \geq U \sqrt{n}$, on a
\begin{eqnarray*}
P(g, g') &=& \Big( \frac{\sin{\phi}}{\phi} \Big)^{\frac{3}{2}}
\frac{\Gamma(\frac{n}{2} + \frac{3}{2})}{\pi^{\frac{n}{2} +
\frac{3}{2}}} 2 \frac{1}{\sqrt{2}} \frac{d_K(g, g')}{D_K(R^2, 2 s)}
B(\frac{n}{2} + 1, \frac{1}{2}) d_K^{-Q-1}(g, g') \\
& &\mbox{} \times \Big( 1 + O(\frac{n}{d_K^2(g, g')}) \Big) \Big( 1
+ O(n^{-\frac{1}{2}}) \Big),
\end{eqnarray*}
o\`u $0 \leq \phi \leq \frac{\pi}{2}$ est d\'etermin\'e par $e^{-
i \phi} = D_K^{-2}(R^2, 2 s) [R^2 - i (2 s)]$.
\end{pro}

{\bf Preuve.} Rappelons que
\begin{eqnarray*}
e^{-\sqrt{-\Delta}} = \frac{1}{2 \sqrt{\pi}} \int_0^{+\infty}
h^{-\frac{3}{2}} e^{-\frac{1}{4h}} e^{h \Delta} \, dh,
\end{eqnarray*}
on ins\`ere \eqref{ehk1} dans la formule pr\'ec\'edente, et par le
th\'eor\`eme de Fubini, on obtient
\begin{eqnarray*}
P(g, g') = \frac{1}{2 \sqrt{\pi}} \frac{1}{(4 \pi)^{\frac{n}{2} +
1}} \int_{\R} \Upsilon \Big( \frac{\lambda}{\sinh{\lambda}}
\Big)^{\frac{n}{2}} \, d\lambda,
\end{eqnarray*}
avec
\begin{eqnarray*}
\Upsilon = \int_0^{+\infty} h^{- \frac{n}{2} - \frac{5}{2}}
\exp{\Big\{ \frac{1}{4 h} \Big[ i (2 s \lambda) - R^2 ( \lambda
\coth{\lambda} - a \frac{\lambda}{\sinh{\lambda}}) -1 \Big] \Big\}}
\, dh.
\end{eqnarray*}

En utilisant le changement de variable $\gamma = \frac{1 + R^2 (
\lambda \coth{\lambda} - a \frac{\lambda}{\sinh{\lambda}}) - i
(2 s \lambda)}{4 h}$, on voit que
\begin{eqnarray*}
\Upsilon = \Big[ \frac{4}{1 + R^2 ( \lambda \coth{\lambda} - a
\frac{\lambda}{\sinh{\lambda}}) - i (2 s \lambda)}
\Big]^{\frac{n}{2} + \frac{3}{2}} \Gamma(\frac{n}{2} + \frac{3}{2}),
\end{eqnarray*}
donc,
\begin{eqnarray*}
P(g, g') &=& \frac{\Gamma(\frac{n}{2} +
\frac{3}{2})}{\pi^{\frac{n}{2} + \frac{3}{2}}} \int_{\R} \Big(
\frac{\lambda}{\sinh{\lambda}} \Big)^{\frac{n}{2}} \Big[ 1 + R^2 (
\lambda \coth{\lambda} - a \frac{\lambda}{\sinh{\lambda}}) - i
(2 s \lambda) \Big]^{-\frac{n}{2} - \frac{3}{2}} \, d\lambda \\
&=& \frac{\Gamma(\frac{n}{2} + \frac{3}{2})}{\pi^{\frac{n}{2} +
\frac{3}{2}}} \int_{\R} \Big( \frac{\sinh{\lambda}}{\lambda}
\Big)^{\frac{3}{2}} \Big\{ \frac{\sinh{\lambda}}{\lambda} \Big[ 1 +
R^2 ( \lambda \coth{\lambda} - a \frac{\lambda}{\sinh{\lambda}}) -
i (2 s \lambda) \Big] \Big\}^{- \frac{n}{2} - \frac{3}{2}} \,
d\lambda.
\end{eqnarray*}

Posons
\begin{eqnarray*}
f(R^2, 2 s, a; \lambda) &=& 1 + R^2 ( \lambda \coth{\lambda} - a
\frac{\lambda}{\sinh{\lambda}}) - i (2 s \lambda) \\
&=& 1 + \frac{\lambda}{\sinh{\lambda}} \Big[ D_K^2(R^2, 2 s)
\cosh{(\lambda - i \phi)} - R^2 a \Big],
\end{eqnarray*}
par \eqref{EGF1}, avec $0 \leq \phi \leq \frac{\pi}{2}$ d\'efini par
$e^{- i \phi} = D_K^{-2}(R^2, 2s) [R^2 - i (2s)]$.

On sait bien que $\Re f(R^2, 2 s, a; \lambda) \geq 1$ pour tout $0
\leq \Im \lambda \leq \frac{\pi}{2}$, voir le Lemma 4.3 de
\cite{P07}. Aussi, pour tout $0 \leq \beta_2 \leq \frac{\pi}{2}$ et
tout $\beta_1 \in \R$, on a
\begin{eqnarray*}
\Re \frac{\sinh{(\beta_1 + i \beta_2)}}{\beta_1 + i
\beta_2} = \frac{\beta_1 \sinh{\beta_1} \cos{\beta_2} + \beta_2
\sin{\beta_2} \cosh{\beta_1}}{\beta_1^2 + \beta_2^2} > 0.
\end{eqnarray*}

D\'efinissons
\begin{eqnarray*}
F(\lambda) = \Big( \frac{\sinh{\lambda}}{\lambda}
\Big)^{\frac{3}{2}} \Big[ \frac{\sinh{\lambda}}{\lambda} +
(D_K^2(R^2, 2 s) \cosh{(\lambda - i \phi)} - R^2 a) \Big]^{-
\frac{n}{2} - \frac{3}{2}},
\end{eqnarray*}
en choisissons la branche principale de la fonction racine car\'ee.
Alors, $F$ est analytique sur
\begin{eqnarray*}
\Omega = \{\lambda \in \C; 0 < \Im \lambda < \frac{\pi}{2} \},
\end{eqnarray*}
et continue sur $\overline{\Omega}$. De plus, on a
\begin{eqnarray*}
\lim_{\lambda \in \overline{\Omega}, |\lambda| \longrightarrow
+\infty} |F(\lambda)| = 0.
\end{eqnarray*}

Par le th\'eor\`eme fondamental de Cauchy, on a donc
\begin{eqnarray*}
\int_{\R} F(\lambda) \, d\lambda = \int_{\R} F(\lambda + i
\phi) \, d\lambda = W,
\end{eqnarray*}
avec
\begin{eqnarray*}
W &=& \int_{\R} \Big( \frac{\sinh{(\lambda + i \phi)}}{\lambda
+ i \phi} \Big)^{\frac{3}{2}} \Big[ \frac{\sinh{(\lambda +
i \phi)}}{\lambda + i \phi} + (D_K^2(R^2, 2 s)
\cosh{\lambda} - R^2 a) \Big]^{- \frac{n}{2} - \frac{3}{2}} \, d\lambda \\
&=& \int_{\R} \Big( \frac{\sinh{(\lambda + i \phi)}}{\lambda +
i \phi} \Big)^{\frac{3}{2}} \Big[ d_K^2(g, g') +
\frac{\sinh{(\lambda + i \phi)}}{\lambda + i \phi} + 2
D_K^2(R^2, 2 s) \sinh^2{\frac{\lambda}{2}} \Big]^{- \frac{n}{2} -
\frac{3}{2}} \, d\lambda.
\end{eqnarray*}

Par cons\'equent, pour terminer la preuve de la Proposition
\ref{PP1}, il nous reste \`a montrer le

\begin{lem} \label{REFL}
Soit $U \gg 1$. Lorsque $n \longrightarrow + \infty$, pour $d_K(g,
g') \geq U \sqrt{n}$, on a
\begin{eqnarray*}
W = \Big( \frac{\sin{\phi}}{\phi} \Big)^{\frac{3}{2}} 2
\frac{1}{\sqrt{2}} \frac{d_K(g, g')}{D_K(R^2, 2 s)} B(\frac{n}{2} +
1, \frac{1}{2}) d_K^{-Q-1}(g, g') \Big( 1 + O(n^{-\frac{1}{2}})
\Big) \Big( 1 + O(\frac{n}{d_K^2(g, g')}) \Big).
\end{eqnarray*}
\end{lem}

{\bf Preuve.} Posons
\begin{eqnarray*}
W_1 &=& \int_{-1}^1 \Big( \frac{\sinh{(\lambda + i
\phi)}}{\lambda + i \phi} \Big)^{\frac{3}{2}} \Big[ d_K^2(g,
g') + \frac{\sinh{(\lambda + i \phi)}}{\lambda + i \phi} +
2 D_K^2(R^2, 2 s) \sinh^2{\frac{\lambda}{2}} \Big]^{- \frac{n}{2} -
\frac{3}{2}} \, d\lambda, \\
W_2 &=& \int_{|\lambda| \geq 1} \Big( \frac{\sinh{(\lambda + i
\phi)}}{\lambda + i \phi} \Big)^{\frac{3}{2}} \Big[ d_K^2(g,
g') + \frac{\sinh{(\lambda + i \phi)}}{\lambda + i \phi} +
2 D_K^2(R^2, 2 s) \sinh^2{\frac{\lambda}{2}} \Big]^{- \frac{n}{2} -
\frac{3}{2}} \, d\lambda.
\end{eqnarray*}

On a alors $W = W_1 + W_2$. On commence par estimer $W_2$. On
constate d'abord que
\begin{eqnarray*}
\Big| \frac{\sinh{(\lambda + i \phi)}}{\lambda + i \phi}
\Big| \leq \cosh{\lambda}, \quad \forall 0 \leq \phi \leq
\frac{\pi}{2}, |\lambda| \geq 1.
\end{eqnarray*}

Donc, pour $n \longrightarrow +\infty$ avec $(\sqrt{2} D_K(R^2, 2 s)
\geq ) d_K(g, g') \gg n^{\frac{1}{2}}$, on a
\begin{eqnarray*}
|W_2| &\leq& \int_{|\lambda| \geq 1} \Big[ (d_K^2(g, g') - 1) + (2
D_K^2(R^2, 2 s) - 2) \sinh^2{\frac{\lambda}{2}} \Big]^{- \frac{n}{2} - \frac{3}{2}} (\cosh{\lambda})^{\frac{3}{2}} \, d\lambda \\
&\leq& 2 \int_1^{+\infty} \Big[ (d_K^2(g, g') - 2) + (2
D_K^2(R^2, 2 s) - 2) \sinh^2{\frac{\lambda}{2}} \Big]^{- \frac{n}{2} - \frac{3}{2}} (\cosh{\lambda})^{\frac{3}{2}} \, d\lambda \\
&\leq& 2 (d_K^2(g, g') - 2)^{- \frac{n}{2} - \frac{3}{2}}
\int_1^{+\infty} (1 + \sinh^2{\frac{\lambda}{2}})^{-\frac{n}{2} - \frac{3}{2}}
(2 \cosh{\frac{\lambda}{2}} \sinh{\frac{\lambda}{2}})^{\frac{3}{2}} \, d\lambda \\
&\leq& (d_K^2(g, g') - 2)^{- \frac{n}{2} - \frac{3}{2}} \frac{16}{n} (\cosh{\frac{1}{2}})^{-n} \\
&=& d_K^{-Q-1}(g, g') e^{-(\frac{n}{2} + \frac{3}{2}) \ln{(1 -
2 d_K^{-2}(g, g'))}} o(n^{-2}) \\
&=& d_K^{-Q-1}(g, g') \Big( 1 + O(\frac{n}{d_K^2(g, g')}) \Big)
o(n^{-2}).
\end{eqnarray*}

Il nous reste \`a \'etudier $W_1$ :

On constate d'abord que pour $-1 \leq \lambda \leq 1$ et $0 \leq
\phi \leq \frac{\pi}{2}$, on a
\begin{eqnarray*}
\Big( \frac{\sinh{(\lambda + i \phi)}}{\lambda + i \phi}
\Big)^{\frac{3}{2}} = \Big( \frac{\sin{\phi}}{\phi}
\Big)^{\frac{3}{2}} \Big[ 1 + O(\lambda) \Big], \qquad
\frac{1}{\cosh{\frac{\lambda}{2}}} = 1 + O(\lambda^2),
\end{eqnarray*}
et si on a de plus $\frac{n}{d_K^2(g, g')} \ll 1$, alors
\begin{eqnarray*}
& &\mbox{} \Big[ d_K^2(g, g') + \frac{\sinh{(\lambda + i
\phi)}}{\lambda + i \phi} + 2 D_K^2(R^2, 2 s)
\sinh^2{\frac{\lambda}{2}} \Big]^{- \frac{n}{2} -
\frac{3}{2}} \\
&=& \Big[ d_K^2(g, g') + 2 D_K^2(R^2, 2 s)
\sinh^2{\frac{\lambda}{2}} \Big]^{- \frac{n}{2} - \frac{3}{2}}
 \\
& & \times \exp\Big\{ -(\frac{n}{2} + \frac{3}{2}) \ln{\Big( 1 +
\frac{\frac{\sinh{(\lambda + i \phi)}}{\lambda + i
\phi}}{1 + \frac{2 D_K^2(R^2, 2 s)}{d_K^2(g, g')}
\sinh^2{\frac{\lambda}{2}}}
\frac{1}{d_K^2(g, g')} \Big)} \Big\} \\
&=& d_K^{-Q-1}(g, g') \Big( 1 + \frac{2 D_K^2(R^2, 2 s)}{d_K^2(g,
g')} \sinh^2{\frac{\lambda}{2}} \Big)^{- \frac{n}{2} - \frac{3}{2}}
\Big[ 1 + O(\frac{n}{d_K^2(g, g')}) \Big] \\
&=& d_K^{-Q-1}(g, g') \Big( 1 + \frac{2 D_K^2(R^2, 2 s)}{d_K^2(g,
g')} \sinh^2{\frac{\lambda}{2}} \Big)^{- \frac{n}{2} - \frac{3}{2}}
\cosh{\frac{\lambda}{2}} \frac{1}{\cosh{\frac{\lambda}{2}}} \Big[ 1
+ O(\frac{n}{d_K^2(g, g')}) \Big].
\end{eqnarray*}

On peut donc \'ecrire
\begin{eqnarray*}
W_1 &=& \Big( \frac{\sin{\phi}}{\phi} \Big)^{\frac{3}{2}}
d_K^{-Q-1}(g, g') \\
& & \times \int_{-1}^1 \Big( 1 + \frac{2 D_K^2(R^2, 2 s)}{d_K^2(g,
g')} \sinh^2{\frac{\lambda}{2}} \Big)^{- \frac{n}{2} - \frac{3}{2}}
\cosh{\frac{\lambda}{2}} \Big[ 1 + O(\frac{n}{d_K^2(g, g')}) \Big]
\Big[ 1 + O(\lambda) \Big] \, d\lambda.
\end{eqnarray*}

Par cons\'equent, pour terminer la preuve du Lemme \ref{REFL}, il
faut montrer que
\begin{eqnarray} \label{EST2}
& &\mbox{} \int_{-1}^1 \Big( 1 + \frac{2 D_K^2(R^2, 2 s)}{d_K^2(g,
g')} \sinh^2{\frac{\lambda}{2}} \Big)^{- \frac{n}{2} - \frac{3}{2}}
\cosh{\frac{\lambda}{2}} \Big[ 1 + O(\frac{n}{d_K^2(g, g')}) \Big]
\Big[ 1 + O(\lambda) \Big] \, d\lambda \nonumber \\
&=& 2 \frac{1}{\sqrt{2}} \frac{d_K(g, g')}{D_K(R^2, 2 s)}
B(\frac{n}{2} + 1, \frac{1}{2}) \Big( 1 + O(n^{-\frac{1}{2}}) \Big)
\Big( 1 + O(\frac{n}{d_K^2(g, g')}) \Big).
\end{eqnarray}

En fait, d'une part, on remarque d'abord que
\begin{eqnarray*}
& &\mbox{} \int_{-1}^1 \Big( 1 + \frac{2 D_K^2(R^2, 2 s)}{d_K^2(g,
g')} \sinh^2{\frac{\lambda}{2}} \Big)^{- \frac{n}{2} - \frac{3}{2}}
\cosh{\frac{\lambda}{2}} \, d\lambda \\
&=& 4 \Big[ \int_0^{+\infty} - \int_{\sinh{\frac{1}{2}}}^{+\infty}
\Big] \Big( 1 + \frac{2 D_K^2(R^2, 2 s)}{d_K^2(g, g')} h^2 \Big)^{-
\frac{n}{2} - \frac{3}{2}} \, dh \\
&=& 4 \frac{1}{\sqrt{2}} \frac{d_K(g, g')}{D_K(R^2, 2 s)} \Big[
\int_0^{+\infty} - \int_{\sqrt{2} \frac{D_K(R^2, 2 s)}{d_K(g, g')}
\sinh{\frac{1}{2}}}^{+\infty} \Big] \Big( 1 +  \lambda^2 \Big)^{-
\frac{n}{2} - \frac{3}{2}} \, d\lambda,
\end{eqnarray*}
puis, par le fait que
\begin{eqnarray*}
\int_0^{+\infty} \Big( 1 +  \lambda^2 \Big)^{- \frac{n}{2} -
\frac{3}{2}} \, d\lambda &=& \frac{1}{2} B(\frac{n}{2} + 1,
\frac{1}{2}), \quad \mbox{(voir \cite{EMOT53} p.10 (16))}, \\
\int_{\sqrt{2} \frac{D_K(R^2, 2 s)}{d_K(g, g')}
\sinh{\frac{1}{2}}}^{+\infty} \Big( 1 +  \lambda^2 \Big)^{-
\frac{n}{2} - \frac{3}{2}} \, d\lambda &\leq& \Big(
\sinh{\frac{1}{2}} \Big)^{-1} \int_{\sinh{\frac{1}{2}}}^{+\infty}
\Big( 1 + \lambda^2 \Big)^{-
\frac{n}{2} - \frac{3}{2}} \lambda \, d\lambda \\
&=& \Big( \sinh{\frac{1}{2}} \Big)^{-1} \frac{1}{n + 1} \Big(
\cosh{\frac{1}{2}} \Big)^{- n - 1},
\end{eqnarray*}
et que
\begin{eqnarray*}
B(\frac{n}{2} + 1, \frac{1}{2}) &=& \Gamma(\frac{1}{2})
\frac{\Gamma(\frac{n}{2} + 1)}{\Gamma(\frac{n}{2} + \frac{3}{2})} =
\sqrt{\pi} \Big( \frac{n}{2} + \frac{3}{2} \Big)^{-\frac{1}{2}} (1
+ o(1)), \\
\mbox{} & & \qquad n \longrightarrow +\infty, \quad \mbox{(voir
\cite{MOS66} p.6 et p.12)},
\end{eqnarray*}
on a
\begin{eqnarray*}
& &\mbox{} \int_{-1}^1 \Big( 1 + \frac{2 D_K^2(R^2, 2 s)}{d_K^2(g,
g')} \sinh^2{\frac{\lambda}{2}} \Big)^{- \frac{n}{2} - \frac{3}{2}}
\cosh{\frac{\lambda}{2}} \, d\lambda \\
&=& 2 \frac{1}{\sqrt{2}} \frac{d_K(g, g')}{D_K(R^2, 2 s)}
B(\frac{n}{2} + 1, \frac{1}{2}) \Big[ 1 + O(n^{-\frac{1}{2}}) \Big],
\quad n \longrightarrow +\infty.
\end{eqnarray*}

D'autre part, puisque $\sinh{h} \geq h$ pour tout $0 \leq h \leq 1$,
on a
\begin{eqnarray*}
& &\mbox{} \int_{-1}^1 |\lambda| \Big( 1 + \frac{2 D_K^2(R^2, 2
s)}{d_K^2(g, g')} \sinh^2{\frac{\lambda}{2}} \Big)^{- \frac{n}{2} -
\frac{3}{2}} \cosh{\frac{\lambda}{2}} \, d\lambda \\
&\leq& 2 \int_0^1 2 \sinh{\frac{\lambda}{2}} \Big( 1 + \frac{2
D_K^2(R^2, 2 s)}{d_K^2(g, g')} \sinh^2{\frac{\lambda}{2}} \Big)^{-
\frac{n}{2} - \frac{3}{2}}
\cosh{\frac{\lambda}{2}} \, d\lambda \\
&\leq& 4 \int_0^{+\infty} \Big( 1 + \frac{2 D_K^2(R^2, 2
s)}{d_K^2(g, g')} y \Big)^{- \frac{n}{2} - \frac{3}{2}} \, dy \leq
\frac{8}{n + 1} \frac{1}{\sqrt{2}} \frac{d_K(g, g')}{D_K(R^2, 2 s)},
\end{eqnarray*}
par le fait que $\frac{1}{\sqrt{2}} \frac{d_K(g, g')}{D_K(R^2, 2 s)}
\leq 1$.

On obtient donc \eqref{EST2}. \cqfd

\medskip

\renewcommand{\theequation}{\thesection.\arabic{equation}}
\section{Preuve du Th\'eor\`eme \ref{TH} pour $M = M_K$ }\label{S4}
\setcounter{equation}{0}

\medskip

Comme
\begin{eqnarray*}
P_h(g) \geq 0, \quad \forall g \in \R^n \times \R, \quad \| P_h \|_1
= 1, \quad \forall h > 0,
\end{eqnarray*}
par ``the Hopf-Dunford-Schwartz maximal ergodic theorem'', on a
\begin{eqnarray*}
\Big| \Big\{g; \sup_{t > 0} \frac{1}{t} \int_0^t e^{-h
\sqrt{-\Delta_G}} f(g) \, dh > \lambda \Big\} \Big| \leq
\frac{2}{\lambda} \| f \|_1, \quad \forall \lambda > 0, f \in
L^1(\R^n \times \R).
\end{eqnarray*}

Pour montrer le Th\'eor\`eme \ref{TH}, il nous reste \`a montrer
qu'il existe une constante $A > 0$ telle que pour $n$ assez grand,
on a
\begin{eqnarray*}
M_K f(g) \leq A n \sup_{t > 0} \frac{1}{t} \int_0^t e^{-h
\sqrt{-\Delta_G}} f(g) \, dh, \quad \forall g \in \R^n \times \R, 0
\leq f \in L^1(\R^n \times \R).
\end{eqnarray*}

Par la structure de dilatation sur $(\R^n \times \R, \Delta_G, dg)$,
il suffit de montrer qu'il existe une constante $L > 0$ telle que
pour $n$ assez grand, avec un $t(n) > 0$ bien choisi, on a pour tout
$g$ et tout $g \neq g' \in B_K(g, 1)$
\begin{eqnarray} \label{EFIN}
| B_K(g, 1) |^{-1} \chi_{B_K(g, 1)}(g') &\leq& L \frac{n}{t(n)}
\int_0^{t(n)} P_h(g, g') \, dh \nonumber \\
&=& L \frac{n}{t(n)} \int_0^{t(n)} h^{-Q} P(\delta_{h^{-1}}(g),
\delta_{h^{-1}}(g')) \, dh.
\end{eqnarray}

Choisissons  $U \gg 1$ et $t(n) = (U \sqrt{n})^{-1}$, pour $0 <
d_K(g, g') < 1$, on constate d'abord que
\begin{eqnarray*}
\frac{1}{t(n)} \int_0^{t(n)} h^{-Q} P(\delta_{h^{-1}}(g),
\delta_{h^{-1}}(g')) \, dh \geq U \sqrt{n} \int_0^{\frac{d_K(g,
g')}{U \sqrt{n}}} h^{-Q} P(\delta_{h^{-1}}(g), \delta_{h^{-1}}(g'))
\, dh,
\end{eqnarray*}
puis, par la Proposition \ref{PP1}, que
\begin{eqnarray*}
& & \frac{1}{t(n)} \int_0^{t(n)} h^{-Q} P(\delta_{h^{-1}}(g),
\delta_{h^{-1}}(g')) \, dh \\
& &\mbox{} \qquad \quad \geq R(n; g, g') U \sqrt{n}
\int_0^{\frac{d_K(g, g')}{U \sqrt{n}}} h \, dh \Big( 1 +
O(\frac{1}{U^2}) \Big) \Big( 1 + O(n^{-\frac{1}{2}}) \Big),
\end{eqnarray*}
avec
\begin{eqnarray*}
R(n; g, g') = \Big( \frac{\sin{\phi}}{\phi} \Big)^{\frac{3}{2}}
\frac{\Gamma(\frac{n}{2} + \frac{3}{2})}{\pi^{\frac{n}{2} +
\frac{3}{2}}} 2 \frac{1}{\sqrt{2}} B(\frac{n}{2} + 1, \frac{1}{2})
\frac{d_K(g, g')}{D_K(R^2, 2 s)} d_K^{-Q-1}(g, g').
\end{eqnarray*}

Si $U$ est choisi assez grand, lorsque $n \longrightarrow +\infty$,
on aura pour tout $0 < d_K(g, g') < 1$,
\begin{eqnarray*}
\frac{1}{t(n)} \int_0^{t(n)} h^{-Q} P(\delta_{h^{-1}}(g),
\delta_{h^{-1}}(g')) \, dh &\geq& \Big( \frac{\sin{\phi}}{\phi}
\Big)^{\frac{3}{2}} \frac{\Gamma(\frac{n}{2} +
\frac{3}{2})}{\pi^{\frac{n}{2} + \frac{3}{2}}} \frac{B(\frac{n}{2} +
1, \frac{1}{2})}{D_K(R^2, 2 s)} \frac{1}{4 U \sqrt{n}} \\
&=& \Big( \frac{\sin{\phi}}{\phi} \Big)^{\frac{3}{2}}
\frac{\Gamma(\frac{n}{2})}{\pi^{\frac{n}{2}}}
\frac{\pi^{-1}}{D_K(R^2, 2 s)} \frac{\sqrt{n}}{8 U},
\end{eqnarray*}
par le fait que $B(\frac{n}{2} + 1, \frac{1}{2}) =
\frac{\Gamma(\frac{1}{2}) \Gamma(\frac{n}{2} +
1)}{\Gamma(\frac{n}{2} + \frac{3}{2})}$ et que $\Gamma(\frac{n}{2} +
1) = \frac{n}{2} \Gamma(\frac{n}{2})$. Or, $|S^{n - 1}| = 2
\frac{\pi^{\frac{n}{2}}}{\Gamma(\frac{n}{2})}$, on a donc
\begin{eqnarray*}
\frac{1}{t(n)} \int_0^{t(n)} h^{-Q} P(\delta_{h^{-1}}(g),
\delta_{h^{-1}}(g')) \, dh \geq \Big( \frac{\sin{\phi}}{\phi}
\Big)^{\frac{3}{2}} \frac{\pi^{-1}}{D_K(R^2, 2 s) |S^{n - 1}|}
\frac{\sqrt{n}}{4 U}.
\end{eqnarray*}

Puisque $\frac{\sin{\phi}}{\phi} \geq \frac{2}{\pi}$ pour tout $0
\leq \phi \leq \frac{\pi}{2}$, pour d\'emontrer \eqref{EFIN}, il
nous reste \`a montrer qu'il existe une constante $C > 0$ telle que
:
\begin{eqnarray} \label{EF1}
n^{\frac{3}{2}} \frac{|B_K(g, 1)|}{D_K(R^2, 2 s) |S^{n - 1}|} \geq
C, \quad \forall g' \in B_K(g, 1), g, n \in \N^*.
\end{eqnarray}

Ce qui est vraie par \eqref{EDK}, \eqref{EVK} et
\begin{eqnarray*}
B(\frac{n}{2}, \frac{3}{2}) \sim n^{-\frac{3}{2}}, \quad n
\longrightarrow +\infty,
\end{eqnarray*}

On a donc montr\'e le Th\'eor\`eme \ref{TH} pour $M = M_K$. \cqfd

\medskip

\renewcommand{\theequation}{\thesection.\arabic{equation}}
\section{Preuve du Th\'eor\`eme \ref{TH} pour $M = M_{CC}$ }\label{S5}
\setcounter{equation}{0}

\medskip

On voit facilement d'apr\`es la section \ref{S4} que pour montrer le
Th\'eor\`eme \ref{TH} pour $M = M_{CC}$, il suffit d'\'etablir les
deux propositions suivantes :

\begin{pro} \label{LL1}
On a $d_K(g, g') \leq d_{CC}(g, g')$ pour tous $g, g'$.
\end{pro}

\begin{pro} \label{LL2}
Il existe une constante $c > 0$ telle que
\begin{eqnarray*}
|B_{CC}(g, 1)| \geq c |B_K(g, 1)|, \quad \forall g, n \in \N^*.
\end{eqnarray*}
\end{pro}

\subsection{Preuve de la Proposition \ref{LL1}}

{\bf I.} Dans le cas o\`u $x' = - x$ et $2 s \geq \pi |x|^2$, on a
\begin{eqnarray*}
d_{CC}^2(g, g') = 2 s \pi, \qquad d_K^2(g, g') = \sqrt{(2 |x|^2)^2 +
(2 s)^2} + 2 |x|^2.
\end{eqnarray*}

On remarque d'abord que $\sqrt{2^2 + \pi^2} \leq \pi^2 - 2$, puis
que $\pi h \geq \sqrt{2^2 + h^2} + 2$ pour tout $h \geq \pi$, on
obtient $d_K(g, g') \leq d_{CC}(g, g')$.

{\bf II.} Dans d'autres cas : Rappelons que
\begin{eqnarray*}
d_K^2(g, g') = \sqrt{R^4 + (2 s)^2} - 2 x' \cdot x, \qquad
d_{CC}^2(g, g') = \Big( \frac{\theta}{\sin{\theta}} \Big)^2 R^2 (1 -
a \cos{\theta}),
\end{eqnarray*}
avec $0 \leq \theta < \pi$ satisfaisant
\begin{eqnarray*}
\frac{2 s}{R^2} = \frac{\theta}{\sin^2{\theta}} - \cot{\theta} + a
\frac{1 - \theta \cot{\theta}}{\sin{\theta}} = \mu(a; \theta) \geq
0.
\end{eqnarray*}

D\'efinissons pour $0 \leq \omega < \pi$, $-1 \leq r \leq 1$,
\begin{eqnarray*}
G(r, \omega) = \Big[ \Big( \frac{\omega}{\sin{\omega}} \Big)^2 (1 -
r \cos{\omega}) + r \Big]^2 - \Big( \frac{\omega}{\sin^2{\omega}} -
\cot{\omega} + r \frac{1 - \omega \cot{\omega}}{\sin{\omega}}
\Big)^2,
\end{eqnarray*}
il nous reste \`a montrer que
\begin{eqnarray} \label{EDCCK}
G(r, \omega) \geq 1, \qquad \forall 0 \leq \omega < \pi, -1 \leq r
\leq 1.
\end{eqnarray}

Pour montrer \eqref{EDCCK}, il suffit d'obtenir les estimations
suivantes :
\begin{eqnarray}
G(-1, \omega) \geq 1, \qquad
\forall 0 \leq \omega < \pi, \label{GAT1}   \\
G(1, \omega) \geq 1, \qquad \forall 0 \leq \omega \leq
\frac{\pi}{2}, \label{GAT2} \\
\inf_{-1 \leq r \leq 1, 0 \leq \omega \leq \frac{\pi}{2}} G(r,
\omega) = \inf_{0 \leq \omega \leq \frac{\pi}{2}} \min\{ G(-1,
\omega), G(1, \omega) \}, \label{GAT3A} \\
\inf_{-1 \leq r \leq 1, \frac{\pi}{2} < \omega < \pi} G(r, \omega) =
\inf_{\frac{\pi}{2} < \omega < \pi} G(-1, \omega). \label{GAT3B}
\end{eqnarray}

La d\'emonstration de \eqref{GAT1} et \eqref{GAT2} est
\'el\'ementaire, mais il y a beaucoup de calcul des d\'eriv\'ees, on
utilisera le logiciel Mathematica pour r\'eduire notre t\^ache.

{\bf Preuve de \eqref{GAT1}.} On constate que
\begin{eqnarray*}
G(-1, \omega) &=& \Big[ \Big( \frac{\omega}{\sin{\omega}} \Big)^2 (1
+ \cos{\omega}) - 1 \Big]^2 - \Big[ \frac{\omega}{\sin^2{\omega}} -
\cot{\omega} - \frac{1 - \omega \cot{\omega}}{\sin{\omega}} \Big]^2
\\
&=& \Big[ 2 \Big( \frac{\frac{\omega}{2}}{\sin{\frac{\omega}{2}}}
\Big)^2 - 1 \Big]^2 - \Big(
\frac{\frac{\omega}{2}}{\sin^2{\frac{\omega}{2}}} -
\cot{\frac{\omega}{2}} \Big)^2.
\end{eqnarray*}

En notant $y = \frac{\omega}{2}$, il nous reste \`a montrer que
\begin{eqnarray*}
Z_1(y) = \Big[ 2 \Big( \frac{y}{\sin{y}} \Big)^2 - 1 \Big]^2 - \Big(
\frac{y}{\sin^2{y}} - \cot{y} \Big)^2 \geq 1, \quad \forall 0 \leq y
< \frac{\pi}{2}.
\end{eqnarray*}

On a \'evidemment $Z_1(0) = 1$, par le th\'eor\`eme des
accroissements finis, il suffit de montrer que $Z_1'(y) \geq 0$ pour
tout $0 < y < \frac{\pi}{2}$. En fait, en utilisant le logiciel
Mathematica, on a
\begin{eqnarray*}
Z_1'(y) = \frac{1}{\sin^5{y}} \Big[ (1 + 6 y^2 - 16 y^4) \cos{y} -
(1 + 2 y^2) \cos{(3 y)} + 2 y (-5 + 8 y^2 + \cos{2 y}) \sin{y}
\Big].
\end{eqnarray*}

Posons
\begin{eqnarray*}
z_*(y) = (1 + 6 y^2 - 16 y^4) \cos{y} - (1 + 2 y^2) \cos{(3 y)} + 2
y (-5 + 8 y^2 + \cos{2 y}) \sin{y}.
\end{eqnarray*}
De la m\^eme fa\c con, comme $z_*(0) = 0$, il nous reste \`a montrer
que $z_*'(y) \geq 0$ pour tout $0 < y < \frac{\pi}{2}$. En utilisant
encore une fois le logiciel Mathematica, on a
\begin{eqnarray*}
z_*'(y) = 4 y \cos{y} (\sin^2{y} - 12 y^2) + 4 y^2 (12 + 4 y^2 + 3
\cos{(2y)}) \sin{y} - 16 \sin^3{y},
\end{eqnarray*}
par un calcul soigneux, on peut \'ecrire (ou bien le v\'erifier par
Mathematica)
\begin{eqnarray*}
z_*'(y) &=& 4 \Big\{ 4 y^2 \sin{y} (y^2 - \sin^2{y}) + \sin{y} (y^2
- \sin^2{y}) \\
& &\mbox{} \quad + (\sin{y} - y \cos{y}) (12 y^2 - 2 y \sin{y}
\cos{y} - 3 \sin^2{y}) \Big\} > 0,
\end{eqnarray*}
pour tout $0 < y < \frac{\pi}{2}$, en utilisant les in\'egalit\'es
\'el\'ementaires
\begin{eqnarray} \label{ECE1}
1 > \frac{\sin{y}}{y} > \Big( \frac{\sin{y}}{y} \Big)^2 > \cos{y},
\qquad \forall 0 < y < \pi.
\end{eqnarray}

Ceci ach\`eve la preuve de \eqref{GAT1}. \cqfd

\medskip

{\bf Preuve de \eqref{GAT2}.} Observons que
\begin{eqnarray*}
G(1, \omega) &=& \Big[ \Big( \frac{\omega}{\sin{\omega}} \Big)^2 (1
- \cos{\omega}) + 1 \Big]^2 - \Big[ \frac{\omega}{\sin^2{\omega}} -
\cot{\omega} + \frac{1 - \omega \cot{\omega}}{\sin{\omega}} \Big]^2
\\
&=& \Big[ 2 \Big( \frac{\frac{\omega}{2}}{\cos{\frac{\omega}{2}}}
\Big)^2 + 1 \Big]^2 - \Big(
\frac{\frac{\omega}{2}}{\cos^2{\frac{\omega}{2}}} +
\tan{\frac{\omega}{2}} \Big)^2.
\end{eqnarray*}

En notant $y = \frac{\omega}{2}$, il nous reste \`a montrer que
\begin{eqnarray*}
Z_2(y) = \Big[ 2 \Big( \frac{y}{\cos{y}} \Big)^2 + 1 \Big]^2 - \Big(
\frac{y}{\cos^2{y}} + \tan{y} \Big)^2 \geq 1, \quad \forall 0 \leq y
< \frac{\pi}{4}.
\end{eqnarray*}

On voit que $Z_2(0) = 0$, et par Mathematica, que
\begin{eqnarray*}
Z_2'(y) = \frac{2}{\cos^5{y}} (\cos{y} + y \sin{y}) (8 y^3 + 2y
\cos{(2y)} - \sin{(2y)}).
\end{eqnarray*}

Il nous reste \`a montrer que
\begin{eqnarray*}
V(h) = h^3 + h \cos{h} - \sin{h} > 0, \qquad \forall 0 < h <
\frac{\pi}{2},
\end{eqnarray*}
ce qui est vraie puisque $V(0) = 0$ et $V'(h) = h( 3 h - \sin{h}) >
0$. \cqfd

\medskip

{\bf Preuve de \eqref{GAT3A}-\eqref{GAT3B}.} D\'efinissons pour $-1
\leq r \leq 1$, $0 \leq \omega < \pi$,
\begin{eqnarray*}
\Phi(r, \omega) = \Big( \frac{\omega}{\sin{\omega}} \Big)^2 (1 - r
\cos{\omega}) + r, \quad \mu(r; \omega) =
\frac{\omega}{\sin^2{\omega}} - \cot{\omega} + r \frac{1 - \omega
\cot{\omega}}{\sin{\omega}}, \\
G_1(r, \omega) = \Phi(r, \omega) - \mu(r; \omega), \qquad G_2(r,
\omega) = \Phi(r, \omega) + \mu(r; \omega).
\end{eqnarray*}

Rappelons que $\mu \geq 0$. On observe d'abord que $\Phi \geq 1$,
puisque pour tout $0 \leq \omega < \pi$,
\begin{eqnarray*}
\frac{\partial \Phi(r, \omega)}{\partial r} = 1 - \Big(
\frac{\omega}{\sin{\omega}} \Big)^2 \cos{\omega} \geq 0, \quad
\Phi(-1, \omega) = 2 \Big(
\frac{\frac{\omega}{2}}{\sin{\frac{\omega}{2}}} \Big)^2 - 1 \geq 1.
\end{eqnarray*}
On a donc $G_2 \geq 1$.

Remarquons que $G = G_1 \cdot G_2$. On va montrer les estimations
suivantes :
\begin{eqnarray} \label{PGAT3}
\frac{\partial G_2(r, \omega)}{\partial r} > 0, \forall 0 < \omega <
\pi, \frac{\partial G_1(r, \omega)}{\partial r} \  \mbox{$\geq 0$
pour $\frac{\pi}{2} \leq \omega < \pi$ et $\leq 0$ pour $0 \leq
\omega \leq \frac{\pi}{2}$.}
\end{eqnarray}

On admet \eqref{PGAT3} pour l'instant et on continue avec la
d\'emonstration de \eqref{GAT3A}-\eqref{GAT3B} :

Pour $0 \leq \omega \leq \frac{\pi}{2}$ fix\'e. Puisque
$\frac{\partial^2 G(r, \omega)}{\partial r^2} = 2 \frac{\partial
G_1(r, \omega)}{\partial r} \frac{\partial G_2(r, \omega)}{\partial
r} \leq 0$, $G(\cdot, \omega)$ est une fonction concave et $\inf_{-1
\leq r \leq 1} G(\cdot, \omega) = \min\{G(-1, \omega), G(1,
\omega)\}$. On obtient donc \eqref{GAT3A}.

Pour $\frac{\pi}{2} < \omega < \pi$ fix\'e. On d\'eduit de
\eqref{PGAT3} que $G_1(r, \omega) \geq G_1(-1, \omega)$. Par le fait
que $G_2 \geq 1$, $G(-1, \omega) = G_1(-1, \omega) G_1(-1, \omega)$
et \eqref{GAT1} impliquent que $G_1(-1, \omega) > 0$ et donc $G_1 >
0$. Par cons\'equent,
\begin{eqnarray*}
\frac{\partial G(r, \omega)}{\partial r} = \frac{\partial G_1(r,
\omega)}{\partial r} G_2(r, \omega) + G_1(r, \omega) \frac{\partial
G_2(r, \omega)}{\partial r} \geq 0,
\end{eqnarray*}
et on obtient \eqref{GAT3B}.

{\bf Preuve de \eqref{PGAT3}.} Observons que
\begin{eqnarray*}
\frac{\partial G_2(r, \omega)}{\partial r} = \Big( 1 -
\frac{\omega^2}{\sin^2{\omega}} \cos{\omega} \Big) +
\frac{1}{\sin{\omega}} \Big( 1 - \frac{\omega}{\sin{\omega}}
\cos{\omega} \Big) > 0,
\end{eqnarray*}
par \eqref{ECE1}.

Remarquons que
\begin{eqnarray*}
\frac{\partial G_1(r, \omega)}{\partial r} = 1 -
\frac{1}{\sin{\omega}} + \frac{\omega}{\sin{\omega}} \cos{\omega}
(\frac{1}{\sin{\omega}} - \frac{\omega}{\sin{\omega}}).
\end{eqnarray*}

Pour $1 < \omega \leq \frac{\pi}{2}$, on a \'evidemment
$\frac{\partial G_1(r, \omega)}{\partial r} \leq 0$.

Pour $0 \leq \omega \leq 1$, par \eqref{ECE1}, on a $\frac{\partial
G_1(r, \omega)}{\partial r} \leq 0$.

Pour $\frac{\pi}{2} \leq \omega < \pi$, pour d\'emontrer
$\frac{\partial G_1(r, \omega)}{\partial r} \geq 0$, en posant
\begin{eqnarray*}
K(\omega) = \sin^2{\omega} - \sin{\omega} - \cos{\omega} (\omega^2 -
\omega),
\end{eqnarray*}
il suffit de montrer que $K(\omega) > 0$ pour tout $\frac{\pi}{2} <
\omega < \pi$. Ce qui est vraie puisque $K(\frac{\pi}{2}) = 0$ et
\begin{eqnarray*}
K'(\omega) = 2 (\sin{\omega} - \omega) \cos{\omega} + (\omega^2 -
\omega) \sin{\omega} > 0, \quad \forall \frac{\pi}{2} < \omega <
\pi.
\end{eqnarray*}

Ceci ach\`eve la preuve de \eqref{PGAT3}. \cqfd

\medskip

\subsection{Preuve de la Proposition \ref{LL2}}

\medskip

On voit que pour tout $u \in \R$ et tout $x \in \R^n$, on a
\begin{eqnarray*}
| B_{CC}((x, u), 1) | = | B_{CC}((x, 0), 1) |, \qquad | B_K((x, u),
1) | = | B_K((x, 0), 1) |.
\end{eqnarray*}

Il suffit de montrer qu'il existe une constante $c > 0$,
ind\'ependante de $n$, telle que
\begin{eqnarray*}
| B_{CC}((x, 0), 1) | \geq c | B_K((x, 0), 1) |, \qquad \forall x
\in \R^n.
\end{eqnarray*}

Lorsque $x = 0$, il suffit de modifier un peu la preuve du Lemme 5.2
de \cite{Li09} dans le cadre des groupes de Heisenberg. On peut donc
supposer dans la suite $x \neq 0$.

D\'efinissons pour $-\pi < \omega < \pi$ et $-1 \leq r \leq 1$,
\begin{eqnarray*}
\Psi(r, \omega) = \Big( \frac{\omega}{\sin{\omega}} \Big)^2 (1 - r
\cos{\omega}), \qquad \Xi(\omega) = \frac{2 \omega - \sin{2
\omega}}{2 \omega^2 \sin{\omega}}.
\end{eqnarray*}

On aura besoin des lemmes suivants :

\begin{lem} \label{DKCL1}
Soit $-1 \leq r \leq 1$. La fonction paire $\Psi(r, \cdot)$ est
strictement croissante sur $[0, \pi)$.
\end{lem}

\begin{lem} \label{DKCL2}
On a {\em\begin{eqnarray} \label{EIF1} \inf_{0 \leq \omega < \pi}
\Xi(\omega) = \Xi(\frac{\pi}{2}) = \frac{2}{\pi}.
\end{eqnarray}}
\end{lem}

\begin{lem} \label{DKCL3}
A l'exception d'un ensemble de mesure nulle, on a {\em
\begin{eqnarray} \label{EIF2}
& &\mbox{} B_{CC}((x, 0), 1) \nonumber\\
&=& \Big\{ (x', t'); \Big( \frac{\omega}{\sin{\omega}} \Big)^2 R^2
(1 - a \cos{\omega}) < 1, t' = \frac{1}{2} \mu(a; \omega) R^2, -\pi
< \omega < \pi \Big\} = \Sigma_1 \nonumber \\
&=& \Big\{ (x', \omega); |x' - x| < 1, -\theta_0 < \omega < \theta_0
< \pi, \  \mbox{avec} \nonumber \\
& &\mbox{} \qquad \qquad  |x' - x|^2 + 2 x' \cdot x (1 -
\cos{\theta_0}) = \Big( \frac{\sin{\theta_0}}{\theta_0} \Big)^2
\Big\} = \Sigma_2.
\end{eqnarray}}
\end{lem}

On admet les lemmes pour l'instant et on continue avec la preuve de
la Proposition \ref{LL2} :

Par le Lemme \ref{DKCL3}, on a
\begin{eqnarray*}
| B_{CC}((x, 0), 1) | = \int \!\!\!\! \int_{\Sigma_1} \, dx'dt' =
\int \!\!\!\! \int_{\Sigma_2} \frac{R^2}{2} \mu'(a; \omega) \,
d\omega dx' = \int_{| z | < 1} R^2 \mu(a; \theta_0) \, dz,
\end{eqnarray*}
o\`u
\begin{eqnarray*}
R^2 = |x|^2 + |x'|^2 = |z|^2 + 2 (|x|^2 + x \cdot z), \qquad a =
\frac{2 (|x|^2 + x \cdot z)}{|z|^2 + 2 (|x|^2 + x \cdot z)},
\end{eqnarray*}
et $0 \leq \theta_0 = \theta_0(x, z) < \pi$ est d\'efini par
\begin{eqnarray} \label{THETA0}
|z|^2 + 2(|x|^2 + x \cdot z) (1 - \cos{\theta_0}) = \Big(
\frac{\sin{\theta_0}}{\theta_0} \Big)^2.
\end{eqnarray}

Donc,
\begin{eqnarray*}
& &\mbox{}|B_{CC}(o, 1)| \\
&=& \int_{| z | < 1} \Big\{ \Big[ |z|^2 +
2(|x|^2 + x \cdot z) \Big] \Big( \frac{\theta_0}{\sin^2{\theta_0}} -
\cot{\theta_0} \Big) + 2(|x|^2 + x \cdot z) \frac{1 - \theta_0
\cot{\theta_0}}{\sin{\theta_0}} \Big\}
\, dz \\
&=& \int_{| z | < 1} \Big\{ \Big[ \Big(
\frac{\sin{\theta_0}}{\theta_0} \Big)^2 + 2(|x|^2 + x \cdot z)
\cos{\theta_0} \Big] \Big( \frac{\theta_0}{\sin^2{\theta_0}} -
\cot{\theta_0} \Big) + 2(|x|^2 + x \cdot z) \frac{1 - \theta_0
\cot{\theta_0}}{\sin{\theta_0}} \Big\} \, dz \\
&=& \int_{| z | < 1} \Big[ \Big( \frac{\sin{\theta_0}}{\theta_0}
\Big)^2 \Big( \frac{\theta_0}{\sin^2{\theta_0}} - \cot{\theta_0}
\Big) + 2(|x|^2 + x \cdot z) \sin{\theta_0} \Big] \, dz.
\end{eqnarray*}

Rappelons que (voir \eqref{VE2} et \eqref{VE1})
\begin{eqnarray*}
4^{-1} (1 + |x|) \int_{| z | < 1} \sqrt{1 - | z |^2} \, dz \leq
|B_K((x, 0), 1)| \leq 2 (1 + |x|) \int_{| z | < 1} \sqrt{1 - | z
|^2} \, dz.
\end{eqnarray*}

Il nous reste \`a montrer qu'il existe une constante $c > 0$,
ind\'ependante de $(n, x, z)$, telle que
\begin{eqnarray} \label{CRSS}
J = J(x, z) = \frac{2 \theta_0 - \sin{2 \theta_0}}{2 \theta_0^2} +
2(|x|^2 + x \cdot z) \sin{\theta_0} \geq c (1 + |x|) \sqrt{1 - | z
|^2}, \  |z| < 1,
\end{eqnarray}
avec $0 \leq \theta_0 < \pi$ d\'etermin\'e par \eqref{THETA0}.

\medskip

\subsection{Preuve de \eqref{CRSS}}

\medskip

Cas 1. $|x| \geq 2$. On d\'eduit de \eqref{THETA0} que
\begin{eqnarray*}
\Big( \frac{\sin{\theta_0}}{\theta_0} \Big)^2 \geq 2 |x| (|x| - 1)
(1 - \cos{\theta_0}),
\end{eqnarray*}
par \eqref{ECE1}, on a $\theta_0 \leq \frac{\pi}{2}$. Donc, d'une
part,
\begin{eqnarray*}
J \geq 2 (|x|^2 + x \cdot z) \sin{\theta_0} \geq 2 |x| (|x| - 1)
\frac{\sin{\theta_0}}{\theta_0} \theta_0 \geq \frac{2}{\pi} |x|^2
\theta_0 \geq \frac{(1 + |x|)}{\pi} |x| \theta_0;
\end{eqnarray*}
d'autre part, \eqref{THETA0} implique que
\begin{eqnarray*}
1 - |z|^2 &=& 1 - \Big( \frac{\sin{\theta_0}}{\theta_0} \Big)^2 + 2
(|x|^2 + x \cdot z) (1 - \cos{\theta_0}) \\
&=& \frac{\theta_0 - \sin{\theta_0}}{\theta_0} \cdot \frac{\theta_0
+ \sin{\theta_0}}{\theta_0} + 4 |x| (|x| + \frac{x}{|x|} \cdot z)
\sin^2{\frac{\theta_0}{2}} \\
&\leq& 2 \frac{\theta_0 - \sin{\theta_0}}{\theta_0} + |x| (|x| + 1)
\theta_0^2, \qquad \mbox{par \eqref{ECE1}}.
\end{eqnarray*}
Or, en utilisant le d\'eveloppement de Taylor d'ordre $3$ autour de
l'origine pour $\sin{\varphi}$, on a
\begin{eqnarray} \label{CREE1}
+\infty > C_1 = \sup_{0 < \omega < \pi} \frac{\omega -
\sin{\omega}}{\omega^3} > \inf_{0 < \omega < \pi} \frac{\omega -
\sin{\omega}}{\omega^3} = c_1 > 0.
\end{eqnarray}
Donc,
\begin{eqnarray*}
1 - |z|^2 \leq \Big[ 2 C_1 + |x| (|x| + 1) \Big] \theta_0^2 \leq 2
(C_1 + 2) |x|^2 \theta_0^2,
\end{eqnarray*}
et on obtient \eqref{CRSS} dans le cas o\`u $|x| \geq 2$.

Cas 2. $|x| < 2$. Il suffit de montrer que
\begin{eqnarray*}
\frac{2 \theta_0 - \sin{2 \theta_0}}{2 \theta_0^2} + 2(|x|^2 + x
\cdot z) \sin{\theta_0} \geq c \sqrt{1 - \Big(
\frac{\sin{\theta_0}}{\theta_0} \Big)^2 + 4 (|x|^2 + x \cdot z)
\sin^2{\frac{\theta_0}{2}}}.
\end{eqnarray*}

Observons que $|x|^2 + x \cdot z \leq |x|^2 + |x| < 6$ et que
\begin{eqnarray*}
1 - \Big( \frac{\sin{\theta_0}}{\theta_0} \Big)^2 + 4 (|x|^2 + x
\cdot z) \sin^2{\frac{\theta_0}{2}} \leq \frac{\theta_0 -
\sin{\theta_0}}{\theta_0} \cdot \frac{\theta_0 +
\sin{\theta_0}}{\theta_0} + 24 \sin^2{\frac{\theta_0}{2}} \leq (2
C_1 + 6) \theta_0^2.
\end{eqnarray*}

Lorsque $|x|^2 + x \cdot z \geq 0$, on a \'evidemment
\begin{eqnarray*}
\frac{2 \theta_0 - \sin{2 \theta_0}}{2 \theta_0^2} + 2(|x|^2 + x
\cdot z) \sin{\theta_0} \geq \frac{2 \theta_0 - \sin{2 \theta_0}}{2
\theta_0^2} \geq c_* \theta_0 > 0, \qquad \forall 0 < \theta_0 <
\pi.
\end{eqnarray*}
Donc, on a l'estimation cherch\'ee.

Lorsque $|x|^2 + x \cdot z < 0$, on constate d'abord que
\begin{eqnarray*}
0 \leq -(|x|^2 + x \cdot z) = |x| (\frac{-x}{|x|} \cdot z - |x|)
\leq \frac{1}{4},
\end{eqnarray*}
puis, par \eqref{EIF1}, que
\begin{eqnarray*}
\frac{2 \theta_0 - \sin{2 \theta_0}}{2 \theta_0^2} + 2(|x|^2 + x
\cdot z) \sin{\theta_0} &=& \sin{\theta_0} \Big[ \Xi(\theta_0) - 2
(- |x|^2 - x \cdot z) \Big] \\
&\geq& (1 - \frac{\pi}{4}) \sin{\theta_0} \Xi(\theta_0) \geq c
\theta_0.
\end{eqnarray*}

Ceci ach\`eve la preuve de \eqref{CRSS}. \cqfd

\medskip

\subsection{Preuve du Lemme \ref{DKCL1}}

\medskip

Lorsque $0 \leq r \leq 1$, il suffit d'utiliser le fait que les deux
fonctions positives $\frac{\omega}{\sin{\omega}}$, $1 - r
\cos{\omega}$ sont strictement croissantes sur $[0, \pi)$.

Pour $-1 \leq r < 0$, il suffit d'observer que
\begin{eqnarray*}
\Big( \frac{\omega}{\sin{\omega}} \Big)^2 (1 - r \cos{\omega}) = (1
+ r) \Big( \frac{\omega}{\sin{\omega}} \Big)^2 + (-2r) \Big(
\frac{\frac{\omega}{2}}{\sin{\frac{\omega}{2}}} \Big)^2.
\end{eqnarray*}

On a donc le Lemme \ref{DKCL1}. \cqfd

\medskip

\subsection{Preuve du Lemme \ref{DKCL2}}

\medskip

Il suffit de montrer que $\Xi'(\omega) < 0$ pour $0 < \omega <
\frac{\pi}{2}$, $\Xi'(\frac{\pi}{2}) = 0$, et $\Xi'(\omega) > 0$
pour $\frac{\pi}{2} < \omega < \pi$. En fait, par un calcul simple,
on a
\begin{eqnarray*}
\Xi'(\omega) = -\frac{\cos{\omega}}{\omega^3 \sin^2{\omega}} \Big[
\omega^2 + \frac{\omega}{2} \sin{(2 \omega)} - 2 \sin^2{\omega}
\Big].
\end{eqnarray*}

Notons
\begin{eqnarray*}
T(\omega) = \omega^2 + \frac{\omega}{2} \sin{(2 \omega)} - 2
\sin^2{\omega}, \quad 0 < \omega < \pi,
\end{eqnarray*}
il nous reste \`a montrer que $T(\omega) > 0$. Ce qui est vraie
puisqu'on a
\begin{eqnarray*}
T(0) = 0, \quad T'(\omega) = 2 \omega + \omega \cos{(2 \omega)} -
\frac{3}{2} \sin{(2 \omega)}, \\
T'(0) = 0, \quad T^{''}(\omega) = 4 \sin{\omega} (\sin{\omega} -
\omega \cos{\omega}) > 0, \quad \forall 0 < \omega < \pi,
\end{eqnarray*}
en utilisant \eqref{ECE1}.  \cqfd

\medskip

\subsection{Preuve du Lemme \ref{DKCL3}}

\medskip

Par le Lemme \ref{DKCL1} et le fait que $\mu(a; \cdot)$ est
strictement croissante, il nous reste \`a montrer que $|z = x' - x|
< 1$. En fait, dans le cas o\`u $|z| < |x|$, il suffit d'observer
que
\begin{eqnarray*}
1 > \Big( \frac{\sin{\omega}}{\omega} \Big)^2 > |x' - x|^2 + 2 x'
\cdot x (1 - \cos{\omega}) \geq |z|^2 + 2 |x| (|x| - |z|) (1 -
\cos{\omega}) \geq |z|^2.
\end{eqnarray*}

Dans le cas o\`u $|z| \geq |x|$, on constate d'abord que
\begin{eqnarray*}
|x' - x|^2 + 2 x' \cdot x (1 - \cos{\omega}) &\geq& |z|^2 \Big[ 1
- 2 \frac{|x|}{|z|} ( 1 - \frac{|x|}{|z|}) (1 - \cos{\omega}) \Big] \\
&\geq& |z|^2 \Big[ 1 - \frac{1}{2} (1 - \cos{\omega}) \Big] = |z|^2
\cos^2{\frac{\omega}{2}},
\end{eqnarray*}
puis que
\begin{eqnarray*}
1 > \Big( \frac{\omega}{\sin{\omega}} \Big)^2 \Big[ |x' - x|^2 + 2
x' \cdot x (1 - \cos{\omega}) \Big] \geq |z|^2 \Big(
\frac{\frac{\omega}{2}}{\sin{\frac{\omega}{2}}} \Big)^2 \geq |z|^2.
\end{eqnarray*}

On obtient donc le Lemme \ref{DKCL3}. \cqfd

\medskip

\renewcommand{\theequation}{\thesection.\arabic{equation}}
\section{Appendice I. Sur l'estimation \eqref{ITI}}
\setcounter{equation}{0}

\medskip

Dans les trois cas connus : les espaces euclidiens $\R^n$, les
groupes de Heisenberg $\mathbb{H}(2 n, 1)$ et les op\'erateurs de
Grushin $\Delta_G$, on peut v\'erifier que pour certaine $U > 0$, il existe
une constante $C(U) > 0$ telle que pour tout $n$ assez grand, on a
\begin{align} \label{ITIA}
\int_0^{\frac{1}{U \sqrt{n}}} e^{-h \sqrt{-\Delta}}(g, g') \, dh
\geq C(U) (- \Delta)^{-\frac{1}{2}}(g, g'), \qquad \forall 0 <
d_K(g, g') < 1.
\end{align}

D\`es qu'on a l'estimation pr\'ec\'edente, on obtient l'estimation
\eqref{ITI} avec une certaine fonction $\phi(n)$. Par la structure
de dilatation et ``the Hopf-Dunford-Schwartz maximal ergodic
theorem'', on a donc
\begin{align} \label{ITIB}
\| M_K \|_{L^1 \longrightarrow L^{1, \infty}} \leq C \phi(n).
\end{align}

Dans \cite{Li09} et aussi dans cet article, l'id\'ee de la
d\'emonstration de \eqref{ITIA} est tr\`es naturelle : il suffit
d'obtenir une estimation inf\'erieure uniform\'ement du noyau de
Poisson. Cependant, ce n'est pas n\'ecessaire : comme on a fait dans
la d\'emonstration, l'\'etape crucial se trouve \`a montrer qu'il
existe deux constantes $c, c' > 0$ telle que pour $n$ assez grand,
on a
\begin{align} \label{ITIC}
(-\Delta)^{-\frac{1}{2}}(g, g') - (-\Delta)^{-\frac{1}{2}} e^{- c
\frac{d_K(g, g')}{\sqrt{n}} \sqrt{-\Delta}}(g, g') \geq c'
(-\Delta)^{-\frac{1}{2}}(g, g'), \quad \forall g \neq g'.
\end{align}

Pour v\'erifier l'estimation \eqref{ITIC} dans le cadre de $\R^n$,
on peut utiliser l'expression explicite de
$(-\Delta)^{-\frac{1}{2}}$ et de $(-\Delta)^{-\frac{1}{2}} e^{-h
\sqrt{-\Delta}}$ ($h
> 0$). Dans le cadre de $\Delta_G$ (resp. $\mathbb{H}(2 n, 1)$), il suffit
de modifier un peu la preuve de la Proposition \ref{PP1} de cet
article (resp. de \cite{Li09}). Dans \cite{LQ11}, on utilisera les
id\'ees ci-dessus afin d'\'etudier la fonction maximale centr\'ee
dans le cadre des groupes de type Heisenberg.

Ce sera beaucoup plus int\'eressant d'\'etudier l'estimation de type
\eqref{ITIC} dans une situation g\'en\'erale (avec $n$ convenable).
Sous les conditions de la propri\'et\'e du doublement du volume et
des estimations gaussiennes classiques du noyau de la chaleur, ce
serait raisonable de croire que \eqref{ITIC} est satisfaite.

De plus, sans utiliser la structure de dilatation, l'estimation
\eqref{ITI} doit s'interpr\'eter par
\begin{align*}
\inf_{n \geq 3, h > 0, g \neq g' \in B(g, h)} \phi(n)
\frac{\sqrt{n}}{h} |B(g, h)| (-\Delta)^{-\frac{1}{2}}(g, g') > 0.
\end{align*}

Enfin, on insiste encore une fois : par notre connaissance, nous ne
savons pas comment \'etablir une estimation de type \eqref{EGGG1}
dans le cadre des groupes de Heisenberg ou dans le cadre des op\'erateurs
de Grushin.

\medskip

\renewcommand{\theequation}{\thesection.\arabic{equation}}
\section{Appendice II. Les cas $p > 1$}
\setcounter{equation}{0}

\medskip

Pour $p > 1$, on a le r\'esultat suivant :

\begin{theo} \label{TTHH1}
Pour tout $1 < p < +\infty$, il existe une constante $A_p
> 0$ telle qu'on a
{\em \begin{eqnarray} \label{Ep} \| M f \|_p \leq A_p \| f \|_p,
\quad \forall f \in L^p, \quad \forall n \in \N^*,
\end{eqnarray}}
avec $M = M_K$ ou bien $M = M_{CC}$.
\end{theo}

{\bf Remarques.} 1. On rappelle qu'une estimation de type \eqref{Ep}
a \'et\'e obtenue par Stein et Str\"omberg dans le cadre des espaces
euclidiens pour la fonction maximale standard de Hardy-Littlewood
(voir \cite{SS83} ou bien  \cite{S82} et \cite{S85}), et par J.
Zienkiewicz dans le cadre des groupes de Heisenberg pour la fonction
maximale d\'efinie par la distance de Carnot-Carath\'eodory ou bien
par celle de Kor\'anyi (voir \cite{Z05}). R\'ecemment, en utilisant l'id\'ee
principale de cet article, l'auteur a
d\'emontr\'e l'estimation \eqref{Ep} dans le cadre des espaces
hyperboliques r\'eels ou complexes, voir \cite{Li11} pour les
d\'etails. Il existe aussi d'autres r\'esulats (partiels), voir par
exemple \cite{B86A}-\cite{C86}, \cite{M90}, \cite{NS94}, \cite{NT09}
et leurs r\'ef\'erences.

2. Remarquons que c'est impossible d'\'etablir une estimation de type \eqref{Ep}
dans le cas g\'en\'eral. Par exemple, pour tout $n \geq 2$ et tout $1 < p_0 < +\infty$,
consid\'erons $\R^+ \times \R^{n - 1}$ muni de la m\'etrique hyperbolique $d_H$
\begin{align*}
d_H((y, x), (y', x')) = \mathrm{arc} \cosh{\frac{y^2 + y'^2 + |x -
x'|^2}{2 y y'}}, \qquad \forall (y, x), (y', x') \in \R^+ \times
\R^{n-1},
\end{align*}
et de la mesure
\begin{align*}
d\mu_{n, p_o}(y, x) = y^{- \frac{p_o}{2 p_o - 1} (n - 1) - 1} \, dydx,
\end{align*}
avec $dx$ la mesure de Lebesgue sur $\R^{n - 1}$.
On sait bien que $(\R^+ \times \R^{n - 1}, d, d\mu_{n, p_o})$ est
\`a croissance exponentielle du volume. Dans cet espace, $M$ est born\'e sur $L^p$ pour $p > p_0$ mais pas pour $1
\leq p < p_0$, voir \cite{Li04} pour les d\'etails et pour plus d'exemples.

3. L'id\'ee de la preuve du th\'eor\`eme \ref{TTHH1} provient de
\cite{Li11}. Pour expliquer \`a peu pr\`es la d\'emonstration, on
observe que
\begin{align*}
\Delta_G =  \sum_{i = 1}^n \frac{\partial^2}{\partial x_i^2} +
(\sum_{i = 1}^n x_i^2) \frac{\partial^2}{\partial u^2},
\end{align*}
\c{c}a nous conduit d'utiliser la s\'eparation des variables et
d'\'etablir une estimation de type
\begin{align} \label{EEMK1}
M_K f(x, u) \leq c M_{\R^n} \Big( M_{\R} f(\cdot, u) \Big)(x),
\end{align}
o\`u la constante $c > 0$ est ind\'ependante de $(n, f, (x, u))$.

\medskip

{\bf Preuve du Th\'eor\`eme \ref{TTHH1}.} Il nous reste \`a
d\'emontrer \eqref{EEMK1}. D\`es qu'on a cette estimation, par les
Propositions \ref{LL1} et \ref{LL2}, on a aussi $M_{CC}f \leq c'
M_{\R^n} ( M_{\R} f)$; ensuite, il suffit d'utiliser le fait que
(voir \cite{S82}, ou \cite{S85}, \cite{SS83}) :
\begin{align*}
\| M_{\R^n} \|_{L^p \longrightarrow L^p} \leq C_p.
\end{align*}

Dans la suite, on donne {\bf la preuve de \eqref{EEMK1}.}

Par la structure de dilatation, il suffit de d\'emontrer que
\begin{align}
\frac{1}{B_K((x, u), 1)} \int_{B_K((x, u), 1)} |f(x', u')| \, dx'du'
\leq c  M_{\R^n} \Big( M_{\R} f(\cdot, u) \Big)(x).
\end{align}

En effet, par \eqref{DDKK1}, on a
\begin{align*}
\int_{B_K((x, u), 1)} |f(x', u')| \, dx'du' &= \int_{B_{\R^n}(x, 1)}
\Big[ \int_{2 |u - u'| < \sqrt{(1 - |x - x'|^2)(1 + |x - x'|^2)}} |f(x', u')| \, du' \Big] \, dx' \\
&\leq \int_{B_{\R^n}(x, 1)} \sqrt{(1 - |x - x'|^2)(1 + |x - x'|^2)} M_{\R}f(x', u) \, dx' \\
&\leq 2 (1 + |x|) \int_{B_{\R^n}(x, 1)} \sqrt{1 - |x - x'|^2}
M_{\R}f(x', u) \, dx',
\end{align*}
en utilisant \eqref{VE2}.

\eqref{VE1} et \eqref{VE2} impliquent que
\begin{align*}
& \frac{1}{B_K((x, u), 1)} \int_{B_K((x, u), 1)} |f(x', u')| \,
dx'du' \\
&\leq 8 \int_{B_{\R^n}(x, 1)}
\sqrt{1 - |x - x'|^2} M_{\R}f(x', u) \, dx' \Big( \int_{B_{\R^n}(x, 1)} \sqrt{1 - |x - x'|^2} \, dx' \Big)^{-1} \\
&= 8 \Big[ M_{\R}f(\cdot, u) \ast \Phi \Big](x),
\end{align*}
avec
\begin{align*}
\Phi(x) = \sqrt{1 - |x|^2} \chi_{B_{\R^n}(o, 1)}(x) \Big(
\int_{B_{\R^n}(o, 1)} \sqrt{1 - |x'|^2} \, dx' \Big)^{-1},
\end{align*}
qui est une fonction non-n\'egative, radialement d\'ecroissante avec
$\| \Phi \|_1 = 1$. Par la propri\'et\'e \'el\'ementaire de
$M_{\R^n}$, on a donc
\begin{align*}
\frac{1}{B_K((x, u), 1)} \int_{B_K((x, u), 1)} |f(x', u')| \, dx'du'
\leq 8 M_{\R^n} \Big( M_{\R} f(\cdot, u) \Big)(x).
\end{align*}

Ceci ach\`eve la preuve du Th\'eor\`eme \ref{TTHH1}. \cqfd

\medskip

Vraisemblablement, en suivant l'id\'ee ci-dessus,
on peut \'etudier l'estimation de type \eqref{Ep} dans le cas o\`u
\begin{align*}
\Delta_G = \Delta_x + |x|^{2 k -2} \Delta_u, \qquad (x, u) \in \R^n \times \R^m, \qquad \mbox {$k, m, n \in \N^*$}.
\end{align*}
Ce sera vraie au moins pour $(m, k)$ fix\'e.

\section*{Remerciements} L'auteur est partiellement support\'e par le
NSF of China (Grant No. 10871048), NCET-09-0316, ``Fok Ying Tong
Education Foundation (Grant No. 111001)'' et ``The Program for
Professor of Special Appointment (Eastern Scholar) at Shanghai
Institutions of Higher Learning''. Je voudrais remercier aussi Prof.
P. Sj\"ogren pour m'avoir indiqu\'e une faute dans l'estimation de
$|B_K(g, 1)|$.

\bigskip

\mbox{}\\
Hong-Quan Li\\
School of Mathematical Sciences  \\
Fudan University \\
220 Handan Road  \\
Shanghai 200433  \\
People's Republic of China \\
E-Mail: hongquan\_li@fudan.edu.cn \quad ou \quad
hong\_quanli@yahoo.fr

\end{document}